\documentclass[11pt]{article}

\usepackage{latexsym,epsfig}
\usepackage{amssymb,amsmath,amsfonts}
\usepackage{exscale} 
\usepackage{ifthen}
\usepackage{longtable}
\usepackage{subfigure}
\usepackage{caption}
\usepackage{cite}
\usepackage{lscape}
\usepackage{alltt} 
\usepackage{multicol} 

\begin{document}
\protect\pagenumbering{arabic}
\setcounter{page}{1}

 \newcommand{\Zt}{\rm}

\newcommand{\ba}{\begin{array}}
\newcommand{\ea}{\end{array}}
\newcommand{\pot}{{\cal P}}
\newcommand{\curv}{\cal C}
\newcommand{\ddt} {\mbox{$\frac{\partial  }{\partial t}$}}
\newcommand{\hl}{\sf}
\newcommand{\hd}{\sf}

\newcommand{\Ad}{\mbox{\rm Ad}}
\newcommand{\Adsm}{\mbox{{\rm \scriptsize Ad}}}
\newcommand{\ad}{\mbox{\rm ad}}
\newcommand{\adsm}{\mbox{{\rm \scriptsize ad}}}
\newcommand{\diag}{\mbox{\rm Diag}}
\newcommand{\sect}{\mbox{\rm sec}}
\newcommand{\id}{\mbox{\rm id}}
\newcommand{\idsm}{\mbox{{\rm \scriptsize id}}}
\newcommand{\eps}{\varepsilon}

\newcommand{\aL}{\mathfrak{a}}
\newcommand{\bL}{\mathfrak{b}}
\newcommand{\mL}{\mathfrak{m}}
\newcommand{\kL}{\mathfrak{k}}
\newcommand{\gL}{\mathfrak{g}}
\newcommand{\nL}{\mathfrak{n}}
\newcommand{\hL}{\mathfrak{h}}
\newcommand{\pL}{\mathfrak{p}}
\newcommand{\uL}{\mathfrak{u}}
\newcommand{\lL}{\mathfrak{l}}

\newcommand{\kG}{{\tt k}}
\newcommand{\nG}{{\tt n}}

\newcommand{\Cart}{$G=K e^{\overline{\aL^+}} K$}
\newcommand{\Area}{\mbox{Area}}
\newcommand{\Hd}{\mbox{\rm Hd}}
\newcommand{\Hdim}{\mbox{\rm dim}_{\mbox{\rm \scriptsize Hd}}}
\newcommand{\Tr}{\mbox{\rm Tr}}
\newcommand{\bs}{{\cal B}}
\newcommand{\nc}{{\cal N}}
\newcommand{\MM}{{\cal M}}
\newcommand{\Ch}{{\cal C}}
\newcommand{\clCh}{\overline{\cal C}}
\newcommand{\Cnt}{\mbox{\rm C}}

\newcommand{\NN}{\mathbb{N}} \newcommand{\ZZ}{\mathbb{Z}}
\newcommand{\QQ}{\mathbb{Q}} \newcommand{\RR}{\mathbb{R}}
\newcommand{\KK}{\mathbb{K}} \newcommand{\FF}{\mathbb{F}}
\newcommand{\CC}{\mathbb{C}} \newcommand{\EE}{\mathbb{E}}
\newcommand{\XX}{X}
\newcommand{\HH}{I\hspace{-2mm}H}
\newcommand{\norm}{\Vert\hspace{-0.35mm}|}
\newcommand{\Sph}{\mathbb{S}}
\newcommand{\ganz}{\overline{\XX}}
\newcommand{\rand}{\partial\XX}
\newcommand{\prodrand}{\partial\XX_1\times\partial\XX_2} 
\newcommand{\regrand}{\partial\XX^{reg}}
\newcommand{\singrand}{\partial\XX^{sing}}
\newcommand{\Frand}{\partial\XX^F}
\newcommand{\Lim}{L_\Gamma}          
\newcommand{\cLim}{M_\Gamma}          
\newcommand{\Flim}{L_\Gamma^B}
\newcommand{\reglim}{L_\Gamma^{reg}}
\newcommand{\radlim}{L_\Gamma^{rad}}
\newcommand{\raylim}{L_\Gamma^{ray}}
\newcommand{\horinf}{\mbox{Vis}^{\infty}}
\newcommand{\horF}{\mbox{Vis}^B}
\newcommand{\Sml}{\mbox{Small}}
\newcommand{\SmlF}{\mbox{Small}^B}

\newcommand{\ifl}{\qquad\Longleftrightarrow\qquad}
\newcommand{\at}{\!\cdot\!}
\newcommand{\ging}{\gamma\in\Gamma}
\newcommand{\xo}{{o}}
\newcommand{\gamo}{{\gamma\xo}}
\newcommand{\gam}{\gamma}
\newcommand{\gax}{h}
\newcommand{\gxi}{{G\!\cdot\!\xi}}
\newcommand{\bd}{$(b,\Gamma\at\xi)$-densit}
\newcommand{\bt}{$(b,\theta)$-densit}
\newcommand{\cd}{$(\alpha,\Gamma\at\xi)$-density}
\newcommand{\be}{\begin{eqnarray*}}
\newcommand{\ee}{\end{eqnarray*}}

\newcommand{\an}{\ \mbox{and}\ }
\newcommand{\as}{\ \mbox{as}\ }
\newcommand{\diam}{\mbox{diam}}
\newcommand{\is}{\mbox{Isom}}
\newcommand{\Ax}{\mbox{Ax}}
\newcommand{\Fix}{\mbox{Fix}}
\newcommand{\Par}{F}
\newcommand{\Min}{\mbox{Fix}}
\newcommand{\vol}{\mbox{vol}}
\newcommand{\Td}{\mbox{Td}}
\newcommand{\piF}{\pi^B}
\newcommand{\piKM}{\pi^I}

\newcommand{\for}{\ \mbox{for}\ }
\newcommand{\pr}{\mbox{pr}}
\newcommand{\sh}{\mbox{sh}}
\newcommand{\shi}{\mbox{sh}^{\infty}}
\newcommand{\rank}{\mbox{rank}}
\newcommand{\supp}{\mbox{supp}}
\newcommand{\mass}{\mbox{mass}}
\newcommand{\kernel}{\mbox{kernel}}
\newcommand{\st}{\mbox{such}\ \mbox{that}\ }
\newcommand{\Stab}{\mbox{Stab}}
\newcommand{\Root}{\Sigma}
\newcommand{\Cone}{\mbox{C}}
\newcommand{\wrt}{\mbox{with}\ \mbox{respect}\ \mbox{to}\ }
\newcommand{\where}{\ \mbox{where}\ }

\newcommand{\con}{{\sc Consequence}\newline}
\newcommand{\rem}{{\sc Remark}\newline}
\newcommand{\prf}{{\sl Proof.\  }}
\newcommand{\qed}{$\hfill\Box$}

\newenvironment{rmk} {\newline{\sc Remark.\ }}{}  
\newenvironment{rmke} {{\sc Remark.\ }}{}  
\newenvironment{rmks} {{\sc Remarks.\ }}{}  
\newenvironment{nt} {{\sc Notation}}{}  

\newtheorem{satz}{\bf Theorem}

\newtheorem{df}{\sc Definition}[section]
\newtheorem{cor}[df]{\sc Corollary}
\newtheorem{thr}[df]{\bf Theorem}
\newtheorem{lem}[df]{\sc Lemma}
\newtheorem{prp}[df]{\sc Proposition}
\newtheorem{ex}{\sc Example}


\title{\sc Geometry and Dynamics of Discrete Isometry Groups of Higher
  Rank Symmetric Spaces}
\author{\sc Gabriele Link}
\date{25/08/2005}
\maketitle
\begin{abstract} For real hyperbolic spaces,  the dynamics of individual
  isometries and the geometry of the limit set of nonelementary discrete
  isometry groups have been studied in great detail. Most of the results were
  generalised to discrete isometry groups of simply connected Riemannian
  manifolds of pinched negative curvature. For symmetric spaces of higher
  rank, which contain isometrically embedded Euclidean planes, the situation
  becomes far more complicated. This paper is devoted to the study of the
  geometric limit set of ``nonelementary'' discrete isometry groups of higher
  rank symmetric spaces. We obtain the natural generalisations of some 
  well-known results from Kleinian group theory. Our main tool consists in  a
  detailed description of the dynamics of individual isometries. As a
  by-product, we give a new geometric construction of free isometry groups with
  parabolic elements in higher rank symmetric spaces.
\end{abstract}
\vspace{0.2cm}

\section{Introduction}

Let $\XX$ be a globally symmetric space
of noncompact type, $\xo\in\XX$  and
$G=\is^o(\XX)$ the connected component of the identity. We will denote by $\rand$ the geometric boundary of $\XX$ endowed with the cone topology (see \cite[chapter~II]{Ba}). 

The goal of this paper is to give more insight into the dynamics of certain individual
isometries of $\XX$ and describe geometrically the structure of the limit set
$\Lim:=\overline{\Gamma\at\xo}\cap\rand$ of discrete isometry groups
$\Gamma\subset G$. The main difficulties we face in the higher rank case
compared to the situation in manifolds with pinched negative curvature are due
to the more complicated structure of the geometric boundary. In fact, to each
point $\xi\in\rand$ we can associate a unique ``direction'' in a fixed Weyl
chamber of $\XX$ (see section~\ref{compactify} for a precise
definition). If the direction of $\xi$ is in the interior of the Weyl chamber,
we say that $\xi$ belongs to the regular boundary $\regrand\subseteq\rand$. Every point in the $G$-orbit of $\xi$ possesses the same
direction, in particular $G$ does not act transitively on the geometric
boundary if the rank of $\XX$ is greater than one. 

Concerning the dynamics of parabolic isometries, there are only  partial results for example by P.~Eberlein (\cite[chapter~4.1]{E}),  A.~Parreau (\cite[chapter~I.2]{P}) and the author (\cite[chapter~4.4]{L}).  
For axial (sometimes also called loxodromic) isometries, however, we are able to describe precisely the action on the geometric boundary. The same remains true for the  particular kind of ``generic parabolic'' isometries, which behave similarly as parabolic isometries of rank one symmetric spaces. A first application gives a new geometric construction of Schottky groups in higher rank symmetric spaces which contain parabolic isometries. 

We then generalise appropriately the notion of ``nonelementary'' groups
well-known in the context of manifolds of pinched negative curvature to higher
rank symmetric spaces. Our notion is weaker than Zariski density  which is normally used for discrete subgroups of real reductive linear Lie groups (see e.g. \cite{B} and \cite{CG}).  Also, the definition is more natural and easily understandable from a geometrical point of view. 

Unfortunately, the  incomplete picture we have about the dynamics of parabolic
isometries makes it difficult to describe the structure of the limit set of
discrete isometry groups which, in general, always contain parabolics. For
the large class of nonelementary groups, however, we can use a so-called
``approximation argument'' in order to reduce the problem to understanding the
action of sequences of axial isometries. A direct consequence of this
approximation argument is Theorem~\ref{axdens}, which states that the set of attractive fixed points of axial isometries is dense in $\Lim$.  It is also one of the key ingredients in the proof of Theorems~\ref{minclosed} and \ref{Product}, which we state here in a simplified version. Let $P_\Gamma$ denote the set of directions of limit points,   and $P_\Gamma^{reg}\subset P_\Gamma$ the set of directions of regular limit points. 
\begin{satz}
If $\Gamma\subset G$ is a nonelementary discrete group, then 
the limit set $K_\Gamma$, considered as a subset of the Furstenberg boundary, is a minimal closed set under the action of $\Gamma$.
\end{satz}
\begin{satz} If $\Gamma\subset G$ is a nonelementary discrete group, then the regular limit set $\reglim$ splits as a product $K_\Gamma\times
  P_\Gamma^{reg}$. 
\end{satz}
To each axial isometry of a higher rank symmetric space we can associate a so-called ``translation vector'', a notion introduced by A.~Parreau (\cite{P}) which generalises the translation length in rank one spaces. 
Let $\ell_\Gamma$ denote the set of translation vectors of axial isometries in $\Gamma$.  Then we have the following
\begin{satz}
 If $\Gamma\subset G$ is a nonelementary discrete group, then  $P_\Gamma$ is equal to the closure of  $\ell_\Gamma$. 
\end{satz}
Although these results are already known for Zariski dense
subgroups of real reductive linear groups (see e.g. \cite{B}, \cite{CG}), the
advantage of our approach is its purely geometric nature which allows to easily adapt the methods to products of manifolds of pinched negative curvature (compare \cite{DK}). 

The paper is organised as follows: In section~2 we recall some basic facts
about Riemannian symmetric spaces of noncompact type and decompositions of semisimple Lie groups. 
We describe the $G$-orbit structure of the geometric boundary
$\rand$ and  introduce a family of (possibly nonsymmetric) $G$-invariant pseudo
distances on $\XX$ which we will need later on. In section~3 we classify individual isometries and describe their dynamics and action on the geometric boundary. As a corollary we obtain Theorem~\ref{freeconst}, which gives a new  geometric  construction of free groups.
Section~4 is devoted to the study of the limit set. We introduce and describe ``nonelementary'' groups of higher rank symmetric spaces. For these groups, we prove the ``approximation argument'' Proposition~\ref{axialapprox} and describe the dynamics of certain sequences of axial isometries. This leads directly to the proof of Theorems~\ref{minclosed}, \ref{Product} and \ref{Cone}.

\section{Preliminaries}

In this section we recall basic facts about symmetric spaces 
of noncompact type (see also \cite{He}, \cite{BGS}, \cite{E}, \cite{W}) and
introduce some notations for the sequel.

\subsection{Cartan and Iwasawa decomposition}\label{cartiwa}
 
Let $\XX$ be a simply connected symmetric space of noncompact type
with base point $\xo\in\XX$, $G=\is^o(\XX)$, and $K$ the isotropy
subgroup of $\xo$ in $G$. It is well-known that $G$ is a semisimple
Lie group with trivial centre and no compact factors, and $K$ a
maximal compact subgroup of $G$. Denote by $\gL$ and $\kL$ the Lie
algebras of $G$ and $K$. Since $G$ acts transitively on $\XX$, we may
identify $\XX$ with the homogeneous space $G/K$. The geodesic symmetry
in $\xo$ induces a Cartan involution $\theta$ on $\gL$, hence
$\gL=\kL\oplus \pL$, where $\pL\subset\gL$ denotes its $-1$
eigenspace. The Killing form $B$ of $\gL$ induces a scalar product 
\begin{equation}\label{scalarproduct}
 \langle X, Y\rangle := - B(X,\theta Y)\,,\quad X,Y\in\gL 
\end{equation}
on $\gL$.  The tangent space $T_\xo\XX$ of $\XX$ in $\xo$ can be
identified with $\pL$, and the Riemannian exponential map at $\xo$ is
a diffeomorphism of $\pL$ onto $\XX$. The scalar product
$\langle\cdot,\cdot\rangle$ of $\gL$  restricted to $\pL$ therefore
induces an scalar product $\langle\cdot,\cdot\rangle$ on  $T_\xo\XX$
which extends to a $G$-invariant Riemannian metric on $\XX$ with associated distance $d$. With respect to this metric, $\XX$ has nonpositive sectional curvature, and, up to rescaling in each factor, this metric is the original one. 

Let $\aL\subset\pL$ be a maximal abelian subspace. Its dimension $r$ is
called the {\hl rank} of $\XX$. The choice of an open  Weyl chamber
$\aL^+\subset \aL$ determines a Cartan decomposition \Cart, where
$\overline{\aL^+}$ denotes the closure of $\aL^+$. We further put $\aL_1:=\{H\in
\aL\;\big\vert\, \Vert H\Vert :=\sqrt{\langle H,H\rangle}=1\}$.

If $z=k e^H\xo\in\XX$, we call $k\in K$ an {\hl angular projection} and
$H\in\overline{\aL^+}$ the (unique) {\hl
  Cartan projection} of $z$. 
\begin{df}\label{thetavec}
For $x,y\in\XX$ the unique vector $H\in\overline{\aL^+}$ with the property 
$x=g\xo$ and $y=g e^H\xo$ for some $g\in G$
is called the {\hd Cartan vector}  of the ordered pair of points  $(x,y)\in\XX\times\XX$ and will be denoted $H(x,y)$. 
\end{df}
Let $\Sigma$ be the set of roots of the pair $(\gL,\aL)$, and
$\Sigma^+\subset\Sigma$ the set of positive roots determined by the
Weyl chamber $\aL^+$. We denote $\gL_\alpha$ the root space of
$\alpha\in\Sigma$, $\nL^+:=\sum_{\alpha\in\Sigma^+} \gL_\alpha$,
and $N^+$ the Lie group exponential of the nilpotent Lie algebra $\nL^+$. 
The decomposition $G=KAN^+$ is called the Iwasawa
  decomposition associated to the Cartan decomposition \Cart. It induces a diffeomorphism $N^+\times \aL\to \XX,\ (n,H)\mapsto
ne^H\xo$, and we have the formula 
\begin{equation}\label{Naction}
d(n e^H\xo,n' e^{H'}\xo)\ge d( e^H\xo,e^{H'}\xo) \qquad\forall\, n,n'\in N^+\
\forall\, H,H'\in \aL \,.
\end{equation}
%
Let $M$ denote the centraliser of $\aL$ in $K$. The Iwasawa decomposition
induces a natural projection
\be \piKM\;:\quad G &\to & K/M\\
g=kan &\mapsto & kM\,, \ee
which we will need in the sequel.

\subsection{The Bruhat decomposition}\label{bruhat}

Given an Iwasawa decomposition $G=KAN^+$, we consider the closed subgroup
$P=MAN^+\subset G$. The homogeneous space $G/P$ is called the {\hl Furstenberg
  boundary} which is identified with $K/M$ via the bijection 
\be \overline{\kappa}\;:\quad G/P &\to & K/M\\ 
gP&\mapsto & \piKM(g)\,.\ee
The Furstenberg boundary hence has a natural differentiable structure arising
from the Lie group structure of $K$. 
Geometrically it can be described as the set of equivalence classes of
asymptotic Weyl chambers in $\XX$ (see \cite{M}).

Let $M^*$ be the normaliser of $\aL$ in $K$, and $W=M^*/M$ the Weyl group
of the pair $(\gL,\aL)$.  We denote by $w_*\in W$ the unique element \st
$\Ad(m_{w_*})(-\aL^+)=\aL^+$ for any representative $m_{w_*}$ of $w_*$ in
$M_*$, and put $\nL^-:=\Ad(m_{w_*})\nL^+=\sum_{\alpha\in\Sigma^+}\gL_{-\alpha}$.  In the sequel we will also need the {\hl opposition involution} 
\begin{eqnarray}\label{oppinv}
 \iota\;:\quad \aL &\to & \aL\nonumber\\
H &\mapsto & -\Ad(m_{w_*}) H\,. 
\end{eqnarray}

The Bruhat decomposition of $G$ \wrt the minimal parabolic subgroup $P$ is 
the disjoint union 
\begin{equation}\label{Bruhatorbit}
 G=\bigcup_{w\in W} N^+ m_w P=\bigcup_{w\in W} U_{w} m_w P\,, 
\end{equation}
where $m_w$ is an arbitrary representative of $w$ in $M^*$, and the sets
$U_w$ are the Lie group exponentials of the subspaces 
$$\uL_w:=\nL^+\cap \Ad(m_w)\nL^-\subset\nL^+\,.$$  
Then the orbit in the decomposition~(\ref{Bruhatorbit}) corresponding to $w_*\in W$ is parametrised by 
$N^+=U_{w_*}$, and the restriction of the above bijection $\overline{\kappa}$
to $N^+ m_{w_*} P$ defines a map
\be \kappa\;:\quad N^+ &\to & K/M\\
n &\mapsto &\overline{\kappa}(n m_{w_*} P)\,. \ee
Geometrically, this map can be interpreted in the following way: If $n\in N^+$, then $\kappa(n)\in K/M$ is the unique element \st  the
Weyl chamber $\kappa(n) e^{\aL^+}\xo$ is asymptotic to the Weyl
chamber $n e^{-\aL^+}\xo$.
The following property of the map $\kappa $ is well-known:
\begin{prp}(\cite{He}, chapter~IX, Corollary 1.9)\label{opendense}\\
The map $\kappa$ is a diffeomorphism onto an open submanifold of $K/M$ 
whose complement consists of finitely many disjoint manifolds of
strictly lower dimension.
\end{prp}
It follows that the orbit $N^+ m_{w_*} P$ is a dense and open
submanifold of the Furstenberg boundary $G/P$. We will call a
$G$-translate of the set $N^+ m_{w_*} P\subset G/P$ a {\hl big cell} of the Furstenberg boundary.

\subsection{Compactification of $\XX$}\label{compactify} 

The geometric boundary $\rand$ of $\XX$ is the set of equivalence
classes of asymptotic geodesic rays endowed with the cone
topology. This boundary is homeomorphic to the unit tangent space of
an arbitrary point in $\XX$.

We fix a Cartan decomposition \Cart\ and let $\xo\in\XX$ be the unique point
stabilised by $K$.  Then for $k\in K$ and $H\in\overline{\aL^+_1}$, we denote by 
$(k,H)$ the unique class in $\rand$ which contains the geodesic ray
$\sigma(t)=ke^{Ht}\xo$, $\, t>0$. We call $k$ an {\hl angular
projection}, and $H$ the {\hl Cartan projection} of
$(k,H)$. Again, the Cartan projection $H$ of a point $\xi\in\rand$ is unique, whereas its angular projection $k$ is only determined up to right multiplication by an element in the centraliser of $H$ in $K$. 

If $r=\rank(\XX)>1$, we define the regular boundary $\regrand$ as the set of
classes with Cartan projection $H\in\aL^+_1$. If $\rank(\XX)=1$, we use the convention $\regrand=\rand$. Furthermore,  the natural 
projection 
$$ \hspace{1cm}\ba{rcl}\pi^B\,:\qquad\ \regrand&\to& K/M\\
(k,H)&\mapsto & kM\, \ea\, $$
will be important in the sequel. The following lemma relates the cone topology to the topology of $K/M$. It is a corollary of Lemma~2.9 in \cite{L}.
\begin{lem}\label{coneKM}
A sequence $(\xi_n)\subset\regrand$ converges to $\xi=(k,H)\in\regrand$ in the cone topology if and only if $\pi^B(\xi_n)$ converges to $kM$ in $K/M$ and the Cartan projections of $\xi_n$ converge to $H$ in $\aL^+_1$. 
\end{lem}
Hence $\pi^B$  is continuous, and $\rank(\XX)=1$ if and only if $\pi^B$ is a homeomorphism.

Let $\ganz:=\XX\cup\rand$. For $x\in\XX$ and $z\in\ganz\setminus \{x\}$ we denote by
$\sigma_{x,z}$ the unique unit speed geodesic emanating from $x$ which
contains $z$. We say that \Cart\ is a Cartan decomposition \wrt $x\in \XX$ (and
$\eta\in\rand$) if $x$ is the unique point fixed by the maximal compact
subgroup $K\subset G$ (and $\sigma_{x,\eta}(t)\subset e^{\overline{\aL^+}}x$
for all $t>0$). 

The isometry group of $\XX$ has a natural action by homeomorphisms on the
geometric boundary. If $g\in G$ and $\xi=(k,H)\in\rand$,  we have $g\at (k,H) = (k_g, H)$,
where $k_g\in K$ is an angular projection of the unique unit speed ray
emanating from $\xo$ asymptotic to the ray
$g\at\sigma_{\xo,\xi}$. Furthermore, the projection $\pi^B$ induces an action
of $G$ by homeomorphisms on the Furstenberg boundary
$K/M=\pi^B(\regrand)$. More precisely, if \Cart\ is a Cartan
decomposition with associated Iwasawa decomposition $G=KAN^+$, and $\pi^I$
the projection defined at the end of section~\ref{cartiwa},
we have the following
\vspace{0mm}\begin{lem}\label{diagact}
Let $g\in G$ and $\xi=(k,H)\in\rand$ with $k\in K$ and
$H\in\overline{\aL^+_1}$. If $k'\in K$ is \st $\pi^I(gk)=k'M$, then $g\xi=
(k', H)$. 

In particular, if $\xi\in\regrand$, then $g\pi^B(\xi)=\pi^B(g\xi)=k'M$.
\end{lem} 
\prf\  Let $\xo\in\XX$ be the fixed point of $K$ and consider the geodesic
$\sigma:=\sigma_{\xo,\xi}$, i.e. $\sigma(t)= k e^{Ht}\xo$ for $t\in\RR$. We write $gk=k'an$ with $k'\in K$, $a\in A$ and $n\in N^+$. In  order to prove that $g \sigma(t) $ converges to $(k',H)\in\rand$ as $t\to\infty$,  
we let $R>>1$ and $\eps>0$ arbitrary. For $t>R\,$ we denote $\sigma_t$ the
geodesic emanating from $\xo$ passing through $g\sigma(t)$. If
$s_t:=d(\xo,g\sigma(t))$, then by the triangle inequality $|s_t-t|\le
d(\xo,g\xo)$. Using the convexity of the distance function we estimate for $t>R$
\be && d(k' e^{HR}\xo, \sigma_t(R))\ \le\ \frac{R}{s_t}\big( d(k'
e^{Hs_t}\xo,g\sigma(s_t))+d(g\sigma(s_t), \sigma_t(s_t))\big)\\
&&\quad =\ \frac{R}{s_t}\big( d(k' e^{Hs_t}\xo,gk e^{Hs_t}\xo)+d(g\sigma(s_t), g\sigma(t))\big)\\
&&\quad =\  \frac{R}{s_t}\big( d(k' e^{Hs_t}\xo,k'an e^{Hs_t}\xo)+d(\sigma(s_t), \sigma(t))\big)
\le \frac{R}{s_t}\big( d(\xo,an\xo)+d(\xo,g\xo)\big)\,\ee
since $d(e^{Hs}\xo, an e^{Hs}\xo)\le d(\xo, an \xo)$ for all $s>0$. From $s_t\to\infty$ as $t\to \infty$ we conclude
$$ d(k' e^{HR}\xo,   \sigma_t(R))<\eps$$ for $t$ sufficiently large.  Hence $g\xi=(k',H)$. \qed\\

Notice that $G\at\xi=K\at\xi$ for any $\xi\in\rand$. Furthermore, $G$ acts
transitively on $\rand$ if and only if $\rank(\XX)=1$.

A further difficulty in the higher rank setting is the fact 
 that, in general, a pair of boundary points can not be joined by a geodesic. The following sets will therefore play a significant role in the sequel.
\begin{df}\label{visisets}
The {\hd visibility set at infinity} viewed from $\xi\in\rand$ is the set
$$ \horinf(\xi):=\{\eta\in\rand\;|\ \exists\ \mbox{geodesic}\ \sigma\  \st 
\sigma(-\infty)=\xi\,,\,\sigma(\infty)=\eta\}\,.$$
The {\hd Bruhat visibility set} viewed from $\xi\in\regrand$ is the image of
$\horinf(\xi)$ under the projection $\pi^B:\;\regrand \to
K/M$, i.e. 
$$ \horF(\xi)=\pi^B(\horinf(\xi))\,.$$
\end{df}
We remark that if  $\rank(\XX)=1$, then $\horF(\xi)\cong\horinf(\xi)=\rand\setminus
\{\xi\}$ for all $\xi\in\rand$. If $\xi\in\regrand$ is stabilised by the minimal parabolic subgroup $P\subset G$, then $\horF(\xi)$ is exactly the image under the map $\overline\kappa$ of the big cell $N^+m_{w_*}P\subset G/P$ of maximal dimension (see \cite[Corollary~2.15]{L} for a more general result). 
In particular, $\horF(\xi)$ can be identified with the nilpotent Lie group
$N^+$ or an arbitrary orbit $N^+ x$, $x\in\XX$. Moreover, 
all Bruhat visibility sets are open and dense submanifolds of $K/M$ by Proposition~\ref{opendense}.  

Furthermore, if $\xi\in\rand$ is arbitrary, $x\in \XX$, and 
$$
N_\xi:=\{g\in G\;|\, \lim_{t\to\infty} d(g\sigma_{x,\xi}(t),\sigma_{x,\xi}(t))=0\}
$$
denotes the {\hl horospherical subgroup} associated to $\xi$, then $\horinf(\xi)=N_\xi\cdot \sigma_{x,\xi}(-\infty)$.

\subsection{Directional distances}\label{directional}

Let $x, y\in \XX$, $\xi\in\rand$, and $\sigma$ a geodesic ray in the
class of $\xi$. We put 
$$ \bs_\xi(x, y)\,:= \lim_{s\to\infty}\big(
d(x,\sigma(s))-d(y,\sigma(s))\big)\,.$$
This number is independent of the chosen ray $\sigma$, and the
function
\be \bs_\xi(\cdot , y):
\quad \XX &\to & \RR\\
x &\mapsto & \bs_\xi(x, y)\ee
is called the {\hl Busemann function} centred at $\xi$ based at $y$ (see also
\cite[chapter~II]{Ba}). 
Using Buseman functions we introduce an important family of (possibly
nonsymmetric) pseudo distances which we will need in the proof of Theorem~\ref{Cone}.
\begin{df}\label{dirdist} Let $\xi\in\rand$. 
We define the {\hl directional distance} of the ordered pair $(x, y)\in
\XX\times \XX$ \wrt the subset $G\at \xi\subseteq\rand$ by
\begin{eqnarray*}
 \bs_{G\cdot\xi}:\quad\XX\times \XX &\to & \RR\\
 (x, y) &\mapsto &  \bs_{G\cdot\xi} (x, y)\;:=\,
 \sup_{g\in  G} \bs_{g \cdot \xi} (x,y)\,.
\end{eqnarray*} 
\end{df}
Notice that in rank one symmetric spaces $G\at\xi=\rand$, hence 
$\bs_{G\cdot\xi}$ equals the Rie\-mannian distance $d$ for any $\xi\in\rand$. In general, the corresponding estimate for the Buseman functions implies
 $$\bs_{G\cdot\xi}(x,y)\le d(x,y)\qquad\forall\,\xi\in\rand\ \
\forall\,x,\,y\in\XX\,.$$
Furthermore, $\bs_{G\cdot\xi}$ is a (possibly nonsymmetric) 
$G$-invariant pseudo distance on $\XX$ (for a proof see \cite[Proposition 3.7]{L}), and we have
$$ \bs_{G\cdot\xi} (x, y) = d(x,y)\cdot \sup_{g\in  G}\ \cos \angle_x
(y, g\xi)\,. $$
In particular, if  \Cart\ is a Cartan decomposition,
$H_\xi\in\overline{\aL^+_1}$  the Cartan projection of $\xi$, and
$H(x,y)\in\overline{\aL^+}$ the Cartan vector of the ordered pair $(x,y)$
according to Definition~\ref{thetavec}, then
\begin{equation}\label{dirwinkel}
 \bs_{G\cdot\xi} (x, y) 
=\langle H_\xi, H(x,y)\rangle \qquad \forall \,x, y\in\XX\,.
\end{equation}

\section{Individual Isometries}\label{individual}

In this section, we recall the  geometric classification and an algebraic characterisation
of individual isometries. We further describe their fixed point set and dynamical properties when acting on the geometric boundary of $\XX$.

\subsection{Geometric classification of isometries}

If $X$ is a Hadamard manifold,  individual
isometries of  $\XX$ can be classified geometrically by means of the
displacement function (compare \cite[chapter~6]{BGS})
$$\ba{rcl} d_\gam\;:\qquad \XX &\to & \RR\\
 x &\mapsto& d(x,\gam  x)\qquad\quad \mbox{for}\quad  \gam\in\is(\XX)
 \,.\ea$$
We will denote by $l(\gam)~:=\inf_{x\in\XX} d_\gam(x)$ the {\hl
 translation  length} of $\gam$.
\begin{df} 
A nontrivial isometry $\gam$ of $\XX$ is called {\hd elliptic}, if
$\gam$ fixes a point in $\XX$, and $\gam$ is called {\hd axial}, if $d_\gam$ assumes
the infimum in $\XX$ and $l(\gam)>0$.

We call $\gam$ {\hd parabolic}, if $d_\gam$ does  not assume the
infimum. If furthermore $l(\gam)=0$, then $\gam$ is called 
{\hd strictly parabolic}, if $l(\gam)>0$, we call $\gam$ {\hd mixed parabolic}.
\end{df}
For $\gam\in\is(\XX)$ we denote by $\Fix(\gam)$ the set of fixed points of $\gam$ in $\ganz$. The following propositions summarise a few properties of individual isometries
of a Hadamard manifold. The proofs can be found in \cite{Ba}, chapter~II.
\begin{prp} (\cite[Proposition~II.3.2]{Ba})\label{ell}\\
An isometry $\gam\in\is(\XX)\setminus\{\id\}$ is elliptic if and only if
$\gam$ has a bounded orbit.
\end{prp}
\begin{prp} (\cite[Proposition~II.3.3]{Ba})\label{ax}\\
An isometry $\gam\in\is(\XX)\setminus\{\id\}$ is axial if and only if there
exists a unit speed geodesic $\sigma$ and a number $l>0$ \st
$\gam(\sigma(t))=\sigma(t+l)$ for all $t\in\RR$. 
\end{prp}
\begin{prp}\label{isopara} (\cite[Proposition~II.3.4]{Ba})\label{par}\\ 
If $\gam\in\is(\XX)\setminus\{\id\}$ is parabolic, then there exists a point 
$\eta\in \Fix(\gam)\subset \partial\XX$ \st $\bs_\eta(x,\gam x)=0$
for all $x\in\XX$.
\end{prp}

\subsection{The Jordan decomposition}\label{Jordan}

From here on we restrict ourselves to the case of a  globally symmetric space $\XX$ of noncompact type. The choice of 
an Iwasawa decomposition $G=KAN^+$ gives rise
to a natural algebraic characterisation of certain individual isometries.
\begin{df}
An isometry $\gam\in G \setminus \{\id\}$ is called
{\hd elliptic}, if $\gam$ is conjugate to an element in $K$,
{\hd hyperbolic}, if $\gam$ is conjugate to an element in $A$, and
{\hd unipotent}, if $\gam$ is conjugate to an element in $N^+$.
\end{df}
Notice that by Proposition~2.19.18 (1), (2) and (5) of \cite{E}, these definitions
coincide with the definitions of ``elliptic'', ``hyperbolic'' and
``unipotent'' via the image under the adjoint representation of $G$ in $GL(\gL)$ .
The following lemma further relates this algebraic characterisation to the geometric classification of the previous subsection. 
\vspace{0mm}\begin{lem}\label{relisos}
$\gam\in G \setminus \{\id\}$ is conjugate to an element in $K$
if and only if $\gam$ fixes a point in $\XX$. Hyperbolic isometries are
axial, and unipotent isometries are strictly parabolic.
\end{lem}
\prf\  The first assertion  is trivial, because the stabiliser of any point in
$\XX$ is conjugate to $K$, the stabiliser of $\xo\in\XX$. 

If $\gam$ is hyperbolic,  there exists $H\in\aL\setminus\{0\}$ and $g\in G$
\st $\gam=g e^H g^{-1}$. The unit speed geodesic $\sigma$ defined by
$\sigma(t):=g e^{t H/||H||}\xo$  then satisfies  
$$\gam (\sigma(t))= (g e^H g^{-1}) g e^{t H/||H||}\xo = g e^{H+ t H/||H||}\xo=\sigma(||H||+t)\qquad \forall\, t\in\RR\,,$$
hence the claim follows from {\rm Proposition}~\ref{ax}. 

If $\gam$ is unipotent, there exists $n\in N^+\setminus\{\id\}$ and $k\in K$
\st $\gam=k n k^{-1}$. For $H\in\aL^+_1$ we define the geodesic $\sigma$ by
$\sigma(t):=k e^{Ht}\xo$, $t\in\RR$. Then 
\be l(\gam)&=&\inf_{x\in\XX} d(x,\gam x) \le \inf_{t>0} d(\sigma(t),\gam \sigma(t))
=\inf_{t>0} d(k e^{Ht}\xo, k n e^{Ht}\xo)\\
&=&\inf_{t>0} d( \xo, e^{-Ht} n e^{Ht}\xo)=0\,,\ee
because  $e^{-Ht} n e^{Ht}\to \id$ as $t \to\infty$. Furthermore, $d_\gam$
does not assume the infimum in $\XX$, because otherwise $\gam$ would be
elliptic.\qed\\

We have the following remarkable Jordan decomposition of elements in $G$.
\begin{thr} (\cite[Theorem~2.19.24]{E})\label{Jordan}\\
For any element $g\in G$ there exists a
unique triplet $e$, $h$, $u$ in $G$ with the following properties:
 $e$ is elliptic, $h$ is  hyperbolic and $u$ is unipotent, 
 $e$, $h$ and $u$ commute pairwise, and  $g=ehu$.
Moreover, if $g'\in G$ commutes with $g$, then $g'$ commutes with $e$, $h$ and $u$. 
\end{thr}
The triplet $e,h,u$ is called the {\hl Jordan decomposition} of $g$. 

Let \Cart\ be a Cartan decomposition. Due to the rich algebraic structure of symmetric spaces, the translation length of an isometry $\gam\in G$ can be  generalised to a vector in $\overline{\aL^+}$. If
$$ C(\gam)~:=\{H(x,\gam x)\;|\,x\in\XX\}\subseteq\overline{\aL^+}\,,$$
where $H(x,\gam x)$ is the Cartan vector from Definition~\ref{thetavec}, then by
Proposition~V.2.1. in~\cite{P}, the closure of $C(\gam)$ in
$\overline{\aL^+}$ contains a unique segment $L(\gam)$ of minimal
length. This segment is called the {\hl translation vector} of
$\gam$. We further have $\Vert L(\gam)\Vert=l(\gam)$ and  $L(\gam)=L(h)$, if $e$, $h$, $u$ is the Jordan
decomposition of $\gam$.   In particular,
$L(\gam)$ is trivial if and only if $\gam$ is elliptic or strictly parabolic.

\subsection{Isometries with positive translation length}\label{postranslen}

We will see in the remainder of this section that there is an essential
difference between isometries $\gamma\in G$ with $l(\gamma)=0$ and
$l(\gamma)>0$. In the first case, we do not know a priori the accumulation
points of an orbit of the cyclic group $\langle\gamma\rangle$ in $\ganz$. On
the other hand, if $l(\gam)$ is positive, then
Proposition~I.2.3 (2) in \cite{P} implies that for all $x\in\XX$ the limit
of the sequence $(\gam^j x)$ as $j\to\infty$ exists and is
independent of $x$. 
\begin{df}\label{attrpoint}
Let $\gam\in G$ be an isometry with positive translation length. The limit 
$\gam^+:=\lim_{j\to\infty} \gam^j\xo$ is called the {\hd attractive fixed point} of $\gam$.  The {\hd repulsive fixed point} $\gam^-$ of $\gam$ is defined as the attractive fixed point of $\gam^{-1}$.  
\end{df}
It will turn out that the fixed points $\gam^+$ and $\gam^-$ play a
significant role in the study of the dynamics of $\gam$. Furthermore, if $e,
h, u$ is the Jordan decomposition of $\gam$, then Proposition~I.2.3
(2) in \cite{P} and its proof show that $\gam^+=h^+$, $\gam^-=h^-$ and $\gam^-\in\horinf(\gam^+)$. 

We fix a Cartan decomposition \Cart\ \wrt $\xo\in\XX$ and the associated Iwasawa decomposition $G=KAN^+$. Recall that $\piKM$ is the natural projection $G\to K/M$ introduced at the end of section~\ref{cartiwa}, and $\iota$ the opposition
involution defined by~(\ref{oppinv}).  The following lemma
describes the coordinates of the fixed points $\gam^+$ and $\gam^-\in\horinf(\gam^+)$. 
\begin{lem} \label{fixpos}
Let  $\gam\in G$ be an isometry with nontrivial translation vector $L\in\overline{\aL^+}$, and $g\in G$   \st the hyperbolic component $h$ in the Jordan decomposition
of $\gam$ belongs to $g e^{\overline{\aL^+}} g^{-1}$. If  $k$,
$k'\in K$ are \st $\pi^I(g)=kM$, $\pi^I(g m_{w_*}^{-1})=k'M$, then 
$\gam^+=(k,L/||L||)$ and $\gam^-=(k',\iota(L)/||L||)$.
\end{lem}
\prf\  By the remark after Definition~\ref{attrpoint} it suffices to prove the
claim for $\gam$ hyperbolic. We first treat the case  $\gam= e^L$. Then
the geodesic $\sigma$ defined by $\sigma(t):= e^{tL/||L||}\xo\ $ for
$t\in\RR\, $ is invariant under $\gam$, 
and we have $\gam^\pm=\sigma(\pm\infty)$.

Furthermore, $\gam^+=\sigma(\infty)=(\id,L/||L||)$,  and for $t>0$ we have
$$\sigma(-t ||L||)= e^{-Lt}\xo=m_{w_*}^{-1} e^{-\Adsm(m_{w_*})Lt }\xo=m_{w_*}^{-1}  e^{\iota(L)
  t}\xo\,,$$ 
hence $\gam^-=\sigma(-\infty)=(m_{w_*}^{-1}, \iota(L)/||L||)$. 

If $\gam=g e^L g^{-1}$ with $g\in G$, then $\gam^\pm=g\sigma(\pm\infty)$, and the assertion follows from {\rm Lemma}~\ref{diagact}. \qed\\

For $\gam\in G$ with $L(\gam)\neq 0$ we put $\ \Par(\gam):=\{x\in\XX\:|\, \sigma_{x,\gam^+}(-\infty)=\gam^-\}\,,$ i.e. 
$\Par(\gam)$ consists of the union of all parallel
geodesics joining $\gam^+$ to $\gam^-$. By Proposition~2.11.4 of \cite{E},
$\Par(\gam)$ is a complete, totally geodesic submanifold of $\XX$. Moreover, if $e$ denotes the elliptic component in the Jordan decomposition of $\gam$, then $\Par(\gam)\cap\Fix(e)\neq\emptyset$.
 
The following proposition shows that if $\gam^+$ and $\gam^-$ are contained in
the regular boundary, then $\gam$ is axial. Although this fact is probably
well-known, we include a geometric proof for the convenience of the reader.   
\vspace{0mm}\begin{prp}\label{transregax}
Let $\gam\in G$ be an isometry with $L(\gam)\in\aL^+\setminus\{\id\}$. Then $\gam$ is axial.
\end{prp} 
\prf\   Let $e,h,u$ be the Jordan decomposition, and
$\gam^+, \gam^-\in\regrand$ the attractive and repulsive fixed point of
$\gam$. We fix a Cartan decomposition \Cart\ \wrt
$\xo\in \Par(\gam)\cap \Fix(e) $ and $\gam^+\in \regrand$, and the
associated Iwasawa  decomposition $G=KAN^+$. By
Proposition~4.1.5 of \cite{E},  $\gam^+$ and $\gam^-$ are fixed by $e$, $h$
and $u$, in particular $e,u\in MAN^+ =\Stab_G(\gam^+)$. Furthermore, $h=e^{L(\gam)}$ and $\Par(\gam)=A\xo$.

We claim that $e\in M$ and $u=\id$. Write $e=man$ with $m\in M$, $a\in A$ and $n\in N^+$. Since $\xo\in\Min(e)$ we have
$$ 0=d(\xo,e\xo)=d(\xo, man\xo)=d(a^{-1}\xo, n\xo)\stackrel{(\ref{Naction})}{\ge} d(a^{-1}\xo,\xo)\,,$$
which implies $a=\id$ and $n=\id$. Hence $e\in M$.

Next we write $u=man$ with $m\in M$, $a\in A$ and $n\in N^+$, and put
$H:=L(\gam)/||L(\gam)||\in \aL_1^+$. We consider the geodesic
$\sigma=\sigma_{\xo,\gam^+}$ and compute for $t>0$
\begin{eqnarray*}
 d(u \sigma(-t), \sigma(-t)) &=& d(man e^{-Ht}\xo, e^{-Ht}\xo) \ge
 d(mane^{-Ht}\xo, ma e^{-Ht}\xo)\\
&& -d(ma  e^{-Ht}\xo, e^{-Ht}\xo)= d(e^{Ht}n e^{-Ht}\xo,\xo)-d(a\xo,\xo)\,.
\end{eqnarray*}
Now if $n\neq \id$, the righthand side is unbounded as $t\to\infty$. Therefore
$u\sigma(-\infty)=\sigma(-\infty)$ implies $n=\id$, hence  
$u=ma$. This is impossible if $u\neq \id$, because the isometry $ma$ assumes the infimum of its displacement function, but $u$ doesn't.\qed

\subsection{Dynamics of axial isometries}\label{dynamicsaxial}

For axial isometries, we are able to describe the action on the geometric
boundary more precisely. 
If $\gam\in\is(\XX)$ is axial, we call the set  
$$\Ax(\gam):=\{x\in\XX\;|\, d(x,\gam
x)=l(\gam)\}$$ the {\hl axis} of $\gam$. It is invariant under the
action of the cyclic group $\langle \gam \rangle$, closed, convex, and consists of
the union of all geodesics translated by $\gam$ (see  
\cite[Proposition 1.9.2 (2)]{E}). In particular, $\Ax(\gam)=\Par(\gam)$, and
the set of fixed points $\Fix(\gam)$ of $\gam$ in $\ganz$ equals $\partial
\Ax(\gam)$.  

We fix a Cartan decomposition \Cart\  of $G=\is^o(\XX)$ with respect to
$\xo\in\XX$, and choose $x\in\Ax(\gamma)$ arbitrary. Then the 
translation vector of $\gam$ is given by 
$$ L(\gam)=H(x,\gam x)\in\overline{\aL^+}\,,$$
where $H(x,\gam x)$ denotes the Cartan vector from Definition~\ref{thetavec}.  
Furthermore, we have the following 
\vspace{0mm}\begin{prp}\label{axial}(\cite[Proposition~2.19.18~(3) and
  Corollary~2.19.19]{E})\\
An isometry $\gam\in G$ is axial if and only if
$\gam$  is conjugate to $e^{H}  m$, where $H\in\overline{\aL^+}\setminus\{0\}$ and
$m\in\{k\in 
K\;|\, \Ad(k)H=H\}$. Moreover,
$H\in\overline{\aL^+}$ equals  the translation vector $L(\gam)$  of $\gam$.  
\end{prp}
For the sake of simplicity, we will restrict ourselves here to the action of
the following kind of axial isometries. For a more general treatment of
arbitrary axial isometries we refer the reader to \cite{L}, chapter~5.2.
\vspace{0mm}\begin{df}
An isometry $\gam\in G$ with translation vector
$L(\gam)\in\aL^+\setminus\{\id\}$ is called a {\hd regular axial} isometry.
\end{df}
In order to describe the dynamics of regular axial isometries on the geometric
boundary $\rand$ of $\XX$, we introduce an auxiliary distance for the Bruhat visibility sets
$\horF(\xi)$, $\xi\in\regrand$, defined in section~\ref{compactify}.

Let \Cart\ be a Cartan decomposition \wrt $x \in\XX$ and $\xi\in\regrand$
arbitrary, and
$G=KAN^+$ the associated Iwasawa decomposition. Recall that $N^+$ is the
Lie group exponential of $\nL^+=\sum_{\alpha\in\Sigma^+} \gL_\alpha$.  
As mentioned at the end of section~\ref{compactify}, $\horF(\xi)$ can be
identified with the submanifold $N^+ x$ of $\XX$.  Let $B_\alpha$
denote the 
scalar product on $\nL^+$  which equals the scalar product~(\ref{scalarproduct})
on $\gL_\alpha$,  and is
zero on $\gL_\beta$ for $\beta\ne \alpha$. Then the scalar product 
$$ ds_{x,\xi}^2 =\frac12 \sum_{\alpha\in\Sigma^+}  B_\alpha \,$$ on $\nL^+$
extends to an $N^+$--invariant metric for the submanifold $N^+ x$ of $\XX$
 with associated Riemannian distance $d_{x,\xi}$ on $N^+ x\cong\horF(\xi)$.  
We remark that for $y\in\XX$, the distance $d_{y,\xi}$ is equivalent to the distance $d_{x,\xi}$ on $\horinf(\xi)$. Furthermore, since the map $\kappa$  is a diffeomorphism from $N^+$ onto a
dense open subset of $K/M$, it follows that the topology induced by the distance $d_{x,\xi}$ on $\horF(\xi)\subset
K/M$ is equivalent to the original topology on $K/M$. The following lemma
describes how this distance behaves under the action of a regular axial isometry. 
Recall that $\iota$ is the opposition involution defined by~(\ref{oppinv}).
\begin{lem}\label{axmove}
  Let $\gax$ be a regular axial isometry, $x\in\Ax(\gax)$ and $\gax^{+}$, $\gax^-$  the attractive and repulsive fixed
  point of $\gax$. Fix a Cartan decomposition \Cart\ \wrt $x$ and $\gax^+$,
  denote $L\in\aL^+$ the translation vector of $\gax$, 
  and put $\alpha_+:=\min \{\alpha(L/||L||)\;|\, \alpha\in\Sigma^+\}$,
  $\alpha_-:=\min \{\alpha(\iota(L)/||L||)\;|\, \alpha\in\Sigma^+\}$.
  If $\eta^+\in\partial\Ax(\gax)\cap\regrand$ satisfies 
  $\piF(\eta^+)=\piF(\gax^+)$, and $\eta^-:=\sigma_{x,\eta^+}(-\infty)$, then 
\begin{eqnarray*}
\forall\,\xi\in\horinf(\eta^+) &\quad & d_{x,\eta^+}(\gax^{-1}\xi, \eta^-)\le
e^{-\alpha_+ ||L||} d_{x,\eta^+}(\xi,
\eta^-)\,,\\
\forall\,\xi\in\horinf(\eta^-)&\quad & d_{x,\eta^-}(\gax\cdot \xi, \eta^+)\le
e^{-\alpha_- ||L||} d_{x,\eta^-}(\xi, \eta^+)\,.
\end{eqnarray*}
\end{lem} 
\prf\  Let $G=KAN^+$ be the Iwasawa decomposition associated to the given Cartan
decomposition.
Then by Proposition~\ref{axial} there exists $m\in M$ \st
$\gax=e^{L}m$. If $\sigma:=\sigma_{x,\eta^+}$, then $\xi\in \horinf(\eta^+)$ implies the existence of $n\in N^+$ \st $\xi=n\sigma(-\infty)$. 

Let $\eps>0$ and $c:[0,1]\to N^+ x\ $ a curve in the submanifold $N^+ x$ with $c(0)=x$, $c(1)=nx$ and
$$ \int_{0}^1 \Vert \dot c (t)\Vert\;dt < d_{x,\eta^+}(\xi, \eta^-)+\eps\,.$$ 
For $t\in [0,1]$, we write $c(t)=n(t)x$ with $n(t)\in N^+$ and put
$$Z(t):=DL_{n(t)^{-1}}\frac{d}{ds}\big|_{s=t} n(s)\in \nL^+\,.$$
Then, by definition of the metric, $\Vert \dot c(t)\Vert^2 = ds_{x,\eta^+}^2(Z(t),Z(t))$. 

Since $\gax^{-1}$ fixes $\eta^-$, and $\gax^{-1}\xi$ corresponds to the element $ \gax^{-1} n\gax x$ in $N^+ x$, the curve $c_\gax(t):=\gax^{-1} n(t) \gax x$
joins  $x$ to $\gax^{-1} n \gax x$,  hence
$$ d_{x,\eta^+}(\gax^{-1}\xi, \eta^-)\le \int_0^1 \Vert \dot{c_\gax}(t)\Vert\;dt\,.$$
Here $ \Vert \dot{c_\gax}(t)\Vert^2 =  ds_{x,\eta^+}^2(\Ad(\gax^{-1}) Z(t),\Ad(\gax^{-1}) Z(t))$.
Since  $M$ normalises $N^+$ we have \mbox{$\Ad(m)Z(t)\in \nL^+$} for all $t\in [0,1]$. Furthermore, 
$$ds_{x,\eta^+}^2(\Ad(m^{-1}) Z(t),\Ad(m^{-1}) Z(t))=ds_{x,\eta^+}^2( Z(t), Z(t))$$
because the scalar product~(\ref{scalarproduct}) on $\gL$ and hence $B_\alpha$, $\alpha\in\Sigma^+$,  is invariant by $\Ad(K)$. 
We conclude
\be ds_{x,\eta^+}^2(\Ad(\gax^{-1}) Z(t),\Ad(\gax^{-1}) Z(t))&= &ds_{x,\eta^+}^2(\Ad(e^{-L}) Z(t),\Ad(e^{-L}) Z(t))\\
&\le&\max_{\alpha\in\Sigma^+} e^{-2\alpha(L) } ds_{x,\eta^+}^2(Z(t),Z(t))\,.\ee
Putting  $\alpha_+:=\min \{\alpha(L/||L||)\;|\, \alpha\in\Sigma^+\}>0$ we summarise
$$ d_{x,\eta^+}(\gax^{-1}\xi, \eta^-)\le e^{-\alpha_+ ||L||} \int_0^1 \Vert \dot c(t)\Vert\;dt < e^{-\alpha_+ ||L||}\left( d_{x,\eta^+}(\xi, \eta^-)+\eps \right)\,,$$
and the first claim follows as $\eps$ tends to zero.

Concerning the second assertion we remark that $\gax^{-1}$ is regular axial with
translation vector 
$\iota(L)\in\aL^+$. Furthermore, $\Ax(\gax^{-1})=\Ax(\gax)$, hence 
the assertion follows from the first claim. \qed\\

For the following important corollary we fix a Cartan decomposition \Cart.
\begin{cor}\label{dist} 
Let $\gax$ be a regular axial isometry and $\sigma\subset\Ax(\gax)$ a
regular geodesic with $\piF(\sigma(\infty))=\piF( \gax^+)$, $\piF(\sigma(-\infty))=\piF(\gax^-)$. Then 
\begin{eqnarray*} 
\forall\,\xi\in\horinf(\sigma(+\infty))&\quad & \lim_{j\to\infty} \gax^{-j}\xi =\sigma(-\infty)\,,\\
\forall\,\xi\in\horinf(\sigma(-\infty))&\quad & \lim_{j\to\infty} \gax^{j}\xi =\sigma(+\infty)\,.
\end{eqnarray*} 
In particular, for all $\xi\in\regrand$ with $\pi^B(\xi)\in \horF(\gax^-)$ the
sequence $\pi^B(\gax^j\xi)$ converges to $\pi^B(\gax^+)$ in $K/M$ as $j\to\infty$. \end{cor}  

\subsection{Dynamics of strictly parabolic isometries} 

We are finally going to study nonelliptic isometries with zero translation
length.  If $\XX$ is a rank one symmetric space, and $\gamma$ a parabolic isometry
of $\XX$, then $\gamma$ fixes a unique point $\eta\in\rand$. Moreover, for all
$z\in\ganz$ we have $\gamma^j z \to
\eta$ and $\gamma^{-j} z\to\eta$ as $j\to \infty$. Unfortunately, this is far from
being true in higher rank symmetric spaces. We only know that 
if $\gamma$ is parabolic with fixed point $\eta$ as in Proposition~\ref{isopara},
then the set of accumulation points of the cyclic group $\langle \gamma\rangle$ is
contained in the boundary of every horosphere centred at $\eta$. Furthermore,
if $\gamma$ is mixed parabolic  then, as pointed out at the beginning of
section~\ref{postranslen}, we have 
$$\lim_{j\to\infty} \gamma^jx = \gamma^+\quad\mbox{and}\qquad
\lim_{j\to\infty}\gamma^{-j}x =\gamma^-\qquad\ \forall \, x\in\XX\,.$$ 
In this section we are going to describe the dynamics of strictly parabolic isometries on the geometric boundary. 

\begin{prp}\label{unipotdyn}
Let $\gam$ be a strictly parabolic isometry, and $\eta\in\rand$ a fixed point of
$\gam$. Then either $\gam$ fixes a point in $\horinf(\eta)$, or for any compact subset $C\subset\horinf(\eta)$ there exists an integer $N\in\NN$ \st 
$$ \gam^{j}C\cap C= \gam^{- j}C\cap C=\emptyset \qquad\forall\, j\ge N\,.$$
\end{prp}
\prf\  Let $e,h,u$ be the Jordan decomposition of $\gam$. Since $\gamma$ is strictly parabolic, we have $h=\id$. Fix a
Cartan decomposition \Cart\ \wrt
$\xo\in \Fix(e) $ and $\eta\in\rand$, and let  $G=KAN^+$ be the associated  Iwasawa  decomposition. Then $u\in N^+$ and $e$ fixes $\eta$. If
$$N_\eta^+:=\{g\in G\;|\, \lim_{t\to\infty} d(g\sigma_{\xo,\eta}(t),
\sigma_{\xo,\eta}(t))=0\}\subseteq N^+$$ 
is the horospherical subgroup associated to $\eta$,  we have $\horinf(\eta)=N_\eta^+\cdot \sigma_{\xo,\eta}(-\infty)$. 
Furthermore, if $H\in\overline{\aL^+_1}$ is the Cartan projection of $\eta$, its Lie algebra is given by 
$$\nL_\eta^+:=\sum_{\begin{smallmatrix}\alpha\in\Sigma^+\\\alpha(H)>0\end{smallmatrix}} \gL_\alpha\subseteq \nL^+\,,$$
and the exponential map $\exp:\nL^+_\eta \to N^+_\eta$ is a diffeomorphism. 
We endow $\nL_\eta^+$ with the norm $\norm\cdot \norm$ associated to the restriction of the scalar product~(\ref{scalarproduct})  to $\nL_\eta^+$, and for $R>0$ we put
$\nL_\eta^+(R):=\{Z\in\nL^+_\eta\;|\, \norm Z\norm<R\}$. 
Then, given a compact set $C\subset\horinf(\eta)$, there exists $R>0$ \st 
$$ C\subset  \;\exp\big(\nL_\eta^+(R)\big)\cdot \sigma_{\xo,\eta}(-\infty)\,.$$
Let $\xi\in\horinf(\eta)$ arbitrary and $Z\in\nL^+_\eta$ \st $\xi=\exp(Z)\cdot\sigma_{\xo,\eta}(-\infty)$. We write 
$$ u=\exp(\sum_{\alpha\in\Sigma^+} Y_\alpha)\,\quad \mbox{and} \qquad
Z=\sum_{\begin{smallmatrix}\alpha\in\Sigma^+\\\alpha(H)>0\end{smallmatrix}}
Z_\alpha\,\quad\mbox{with}\qquad   Y_\alpha,\, Z_\alpha\in\gL_\alpha\,.$$
If  $Y_\alpha=0$ for all $\alpha\in\Sigma^+$ with $\alpha(H)>0$, then
$\ad(H)Y_\alpha=0$ for all $\alpha\in \Sigma^+$, hence $e^{Ht}ue^{-Ht}= u$ for
all $t\in\RR$. Since $e$ fixes $\sigma_{\xo,\eta}$ pointwise, we conclude
$$ d(\sigma_{\xo,\eta}(-t),\gamma \sigma_{\xo,\eta}(-t))= d(\sigma_{\xo,\eta}(-t),u \sigma_{\xo,\eta}(-t))=d(\xo, e^{H t}u e^{-Ht}\xo)= d(\xo, u\xo)$$
for all $t\in\RR$. In particular, $\gamma$ fixes $\sigma_{\xo,\eta}(-\infty)$.

So if $\gamma$ does not fix a point in $\horinf(\eta)$, there exists
$\alpha\in\Sigma^+$ \st $\alpha(H)>0$ and $Y_\alpha\ne 0$. We may therefore choose $\beta\in \Sigma^+$ \st $\beta(H)>0$, $Y_\beta\neq 0$  and $\beta(H) \le \alpha(H)$ for all $\alpha\in\Sigma^+$ with $\alpha(H)>0$ and 
$Y_\alpha\neq 0$. For $j\in\NN$ we write 
$$u^{\pm j} \exp(Z) = \exp(Y^{(\pm j)})=\exp(\sum_{\alpha\in\Sigma^+} Y_\alpha^{(\pm
  j)})\,\quad\mbox{with}\qquad Y_\alpha^{(\pm j)}\in\gL_\alpha\,.$$ 
Then $[\gL_\alpha,\gL_{\alpha'}]\subseteq \gL_{\alpha+\alpha'}$ for all
  $\alpha, \alpha'\in \Sigma^+$  and the Campbell Hausdorff formula imply 
$$ Y_{\beta}^{(\pm j)}=\pm j Y_\beta+ Z_\beta+ X_\beta\,,$$ 
where $X_\beta\in \gL_\beta$ is a term consisting of successive Lie brackets
of the $Z_\alpha$. In particular, $\norm X_\beta\norm$ is bounded, and
therefore 
$\norm Y_{\beta}^{(j)}\norm$ and $\norm Y_\beta^{(-j)}\norm$ tend to infinity as
$j\to\infty$. Since root spaces associated to different roots
in $\Sigma^+$ are orthogonal \wrt the scalar product~(\ref{scalarproduct}) (see \cite[Theorem~III.4.2 (iii)]{He}), we have  
$$ \norm Y^{(\pm j)}\norm = \sqrt{\sum_{\alpha\in\Sigma^+} \norm Y_\alpha^{(\pm
  j)}\norm^2}  \ge \norm Y_{\beta}^{(\pm j)}\norm \ge  j  \norm Y_\beta\norm -\norm Z_\beta\norm -\norm X_\beta\norm  \,. $$ 
Now by $\Ad(K)$-invariance of the scalar product~(\ref{scalarproduct}) and $e\in K$ we obtain 
$\norm \Ad(e^{\pm j}) Y^{(\pm j)}\norm =\norm  Y^{(\pm j)}\norm \ge \norm  Y_\beta^{(\pm j)}\norm$.  Hence if $j>( R+\norm Z_\beta\norm +\norm X_\beta\norm)/\norm Y_\beta\norm$, we deduce
$$ \gamma^{\pm j} \exp(Z)\sigma_{\xo,\eta}(-\infty)=  \exp(\Ad(e^{\pm j}) Y^{(\pm j)})\sigma_{\xo,\eta}(-\infty)\notin C\,,$$
since $e$ commutes with $u$   and fixes $\sigma_{\xo,\eta}$. By compactness of
$C$ there exists $N\in\NN$ \st for any $\xi\in C$ and all $j\ge N$ we have 
$\ \gamma^{j}\xi\notin C$ and $\ \gamma^{-j}\xi\notin C$.\qed\\

For the following kind of parabolic isometries we have similar dynamics on the
geometric boundary as in rank one symmetric spaces. 
They will also play an important role for the construction of free groups in
the following section. 
\begin{df}
An isometry $\gam\in G \setminus \{\id\}$ is called
{\hd generic parabolic}, if $\gamma$ is strictly parabolic and possesses a unique fixed point in each $G$--invariant subset of the regular geometric boundary $\regrand$.
\end{df}
The following lemma will be convenient in the sequel. 
\begin{lem}
If $\gamma$ is generic parabolic and $e,h,u$ its Jordan decomposition, then $h=\id$, $e$ is conjugate to an element in $M$, and $u$ is conjugate to an element $n\in N^+$ with the property $n\notin m_w N^+m_w^{-1}$ for all nontrivial Weyl group elements
$w\in W\setminus\{\id\}$. 
\end{lem}
\prf\  Let $e,h,u$ be the Jordan decomposition and $\eta\in\regrand$  a fixed point of $\gamma$. Since $\gamma$ is strictly parabolic, we have $l(\gamma)=l(h)=0$, hence $h=\id$. Furthermore, Proposition~4.1.5 of \cite{E} implies that $\eta$ is fixed by both $e$ and $u$. 
Let \Cart\ be a Cartan decomposition \wrt $\xo\in\Min(e)$ and
$\eta\in\regrand$, and $G=KAN^+$ the associated Iwasawa decomposition. Then $u\in N^+$ and $e\in MAN^+=\Stab_G(\eta)$.
We write $e=man$ with $m\in M$, $a\in A$ and $n\in N^+$.  Since $\xo\in\Min(e)$ we have
$$ 0=d(\xo,e\xo)=d(\xo, man\xo)=d(a^{-1}\xo, n\xo)\stackrel{(\ref{Naction})}{\ge} d(a^{-1}\xo,\xo)\,$$
which implies $a=\id$ and $n=\id$. Hence $e\in M$.

Next suppose there exists $w\in W\setminus\{\id\}$ \st $u\in m_w N^+ m_w^{-1}$. Then $m_w\eta\neq \eta$ and, since $N^+$ fixes $\eta$, we have $u m_w\eta =m_w\eta$. Now $e\in M$ implies $e m_w\eta=m_w\eta$,  i.e. $m_w\eta\in G\at\eta$ is fixed by $\gamma=eu$, in contradiction to the fact that $\eta$ is the unique fixed point of $\gamma$ in $G\at\eta$. \qed \\ 

Recall the definition of
the map $\overline{\kappa}$ from section~\ref{bruhat}.
\begin{prp}\label{dyngenuni}
If $\gam$ is generic parabolic and $\eta\in\regrand$ a fixed point of $\gam$, then for any $\xi\in G\at\eta$ we have 
$$ \lim_{j\to\infty} \gam^j\xi=\lim_{j\to\infty}\gam^{-j}\xi=\eta\,.$$
\end{prp}
\prf\ Let $e, h=\id , u$  be the Jordan decomposition of $\gamma$. We fix a
Cartan decomposition \Cart\ \wrt $\xo\in\Min(e)$ and $\eta\in\regrand$, and the
associated Iwasawa decomposition $G=KAN^+$. Denote by $H\in\aL^+_1$ the
Cartan projection of $\eta$ and let $V\subset K/M$ be an open neighbourhood of $\piF(\eta)$. For $w\in W$
we denote $\Vert \cdot \Vert_{\begin{smallmatrix} \\\tiny w\end{smallmatrix}}$ the
norm on $\uL_w:=\nL^+\cap \Ad(m_w )\nL^-$ associated to the restriction of the scalar product~(\ref{scalarproduct}) to $\uL_w$. Then by the Bruhat decomposition there
exists $R>0$ \st with $U_w(R):=\{\exp Z\;|\, Z\in\uL_w\ \mbox{with}\
||Z||_{\begin{smallmatrix} \\\tiny w\end{smallmatrix}}\le R\}$ we have 
$$ K/M\setminus \overline\kappa \big(\bigcup_{w\in W\setminus\{\idsm\}} U_w(R) m_w
P\big)\subset V\,.$$
Hence it suffices to prove that for all $\xi\in \regrand$ and $j$ sufficiently
large
$$\piF(\gam^{\pm j}\xi)\in K/M\setminus \overline\kappa \big(\bigcup_{w\in
  W\setminus\{\idsm\}} U_w(R) m_w P\big)\,.$$
Let $\xi\in\regrand$ arbitrary. If $\piF(\xi)=\piF(\eta)$, there is nothing to prove. If $\piF(\xi)\ne \piF(\eta)$, there exists $w\in W\setminus\{\id\}$ and $Z\in
\uL_w$ \st $\piF(\xi)\in \overline{\kappa} (\exp(Z) m_w P)$.  As in the proof of the
previous proposition we write
$$ u=\exp(\sum_{\alpha\in\Sigma^+} Y_\alpha)\,\quad \mbox{and} \qquad
Z=\sum_{\alpha\in\Sigma^+}
Z_\alpha\,\quad\mbox{with}\qquad   Y_\alpha,\, Z_\alpha\in\gL_\alpha\,.$$
Since $u\notin m_wN^+m_w^{-1}$ there exists $\alpha\in\Sigma^+$ \st
$Y_\alpha\ne 0$ and $\Ad(m_w^{-1})Y_\alpha\in \nL^-$,
i.e. $Y_\alpha\in\uL_w$. Hence we may choose $\beta\in \Sigma^+$ \st $\beta(H)\le \alpha(H)$ for all $\alpha\in\Sigma^+$ with
$Y_\alpha\neq 0$ and $Y_\alpha\in\uL_w$. For $j\in\NN$, we write 
$$u^{\pm j} \exp(Z)= \exp(Y^{(\pm j)})=\exp(\sum_{\alpha\in\Sigma^+} Y_\alpha^{(\pm
  j)})\,\quad\mbox{with}\qquad Y_\alpha^{(\pm j)}\in\gL_\alpha\,.$$ 
Then $[\gL_\alpha,\gL_{\alpha'}]\subseteq \gL_{\alpha+\alpha'}$ for all
  $\alpha, \alpha'\in \Sigma^+$  and the Campbell Hausdorff formula imply 
$$ Y_{\beta}^{(\pm j)}=\pm j Y_\beta+ Z_\beta+ X_\beta\,,$$ 
where $X_\beta\in \gL_\beta$ is a term consisting of successive Lie brackets
of the $Z_\alpha$. In particular, $\Vert X_\beta\Vert_{\begin{smallmatrix} \\\tiny w\end{smallmatrix}}$ is bounded, and
therefore 
$\Vert Y_{\beta}^{(j)}\Vert_{\begin{smallmatrix} \\\tiny w\end{smallmatrix}}$ and $\Vert Y_\beta^{(-j)}\Vert_{\begin{smallmatrix} \\\tiny w\end{smallmatrix}}$ tend to infinity as
$j\to\infty$. As before, we obtain by the $\Ad(K)$-invariance of the scalar product~(\ref{scalarproduct})
$$ \Vert \Ad(e^{\pm j}) Y^{(\pm j)}\Vert_{\begin{smallmatrix} \\\tiny w\end{smallmatrix}}=\Vert Y^{(\pm j)}\Vert_{\begin{smallmatrix} \\\tiny w\end{smallmatrix}} \ge \Vert Y_{\beta}^{(\pm j)}\Vert_{\begin{smallmatrix} \\\tiny w\end{smallmatrix}} \ge j  \Vert Y_{\beta}\Vert_{\begin{smallmatrix} \\\tiny w\end{smallmatrix}}-\Vert Z_{\beta}\Vert_{\begin{smallmatrix} \\\tiny w\end{smallmatrix}}-\Vert X_{\beta}\Vert_{\begin{smallmatrix} \\\tiny w\end{smallmatrix}}\,. $$ 
Hence for $j\in\NN$ sufficiently large we have $ \Vert  \Ad(e^{\pm j}) Y^{(\pm j)}\Vert_{\begin{smallmatrix} \\\tiny w\end{smallmatrix}}>R$ and the claim follows from 
$\piF(\gam^{(\pm j)}\xi)= \overline\kappa\big(\exp( \Ad(e^{\pm j}) Y^{(\pm j)})m_w P\big)$.   \qed\\

\subsection{Construction of free groups}\label{free}
 
We will now apply the results of the previous sections to produce
Schottky groups, an interesting kind of free and discrete isometry groups of infinite covolume. Their construction is based  on the following
\vspace{0mm}\begin{lem}[Klein's Criterion](see \cite{Ha})\\
Let $G$ be a group acting on a set $S$, $\Gamma_1,\Gamma_2$ two subgroups of $G$, where $\Gamma_1$ contains at least three elements, and let $\Gamma$ be the subgroup they generate.  Assume that there exist two nonempty subsets $S_1$, $S_2$ in $S$ with $S_2$ not included in $S_1$  \st 
$\gamma (S_2)\subseteq S_1$ for all $\gamma\in\Gamma_1\setminus\{\id\}$ and
$\gamma (S_1)\subseteq S_2$ for all $\gamma\in\Gamma_2\setminus\{\id\}$.
Then $\Gamma$ is isomorphic to the free product $\Gamma_1 * \Gamma_2$.
 \end{lem}
For the remainder of this section we fix a Cartan decomposition \Cart\ \wrt the base point $\xo\in\XX$.
Recall from the remark following Definition~\ref{visisets} that a finite
intersection of sets $\horF(\xi_i)\subset K/M$, $\xi_i\in\regrand$, is a dense
and open subset of $K/M$. The following theorem  describes a new geometric
construction of finitely generated free groups containing parabolic isometries. 
\begin{thr}\label{freeconst}
Let $\XX$ be a globally symmetric space of noncompact type, and\\
$\{\xi_1,\xi_2,\ldots, \xi_{2l}, \xi_{2l+1},\ldots,
\xi_{2l+p}\}\subset\regrand$ a set of $\,2l+p$ points \st 
$$ \pi^B(\xi_i)  \in  \bigcap_{\begin{smallmatrix}n=1\\
n\neq i\end{smallmatrix}}^{2l+p} \horF(\xi_n)\qquad\forall\, i\in\{1,2,\ldots, 2l+p\}\,.$$
Then there exist regular axial isometries $\gamma_1, \gamma_2,\ldots,\gamma_l$
with $\gamma_m^+=\xi_{2m}$ and
$\piF(\gamma_m^-)=\piF(\xi_{2m-1})\,$ for $1\le m\le l$, and generic parabolic
isometries  $\gamma_{l+1},\gamma_{l+2},\ldots,\gamma_{l+p}$ with fixed points $\xi_{2l+1},
\xi_{2l+2},\ldots,\xi_{2l+p}$ respectively.  
Furthermore, there exist pairwise disjoint open neighbourhoods $U_i\subset
K/M$ of $\piF(\xi_i)$, $1\le i\le 2l+p$, \st with $C:=\bigcup_{i=1}^{2l+p} \overline{U_i}$ we
have 
\be && \gamma_m(C\setminus \overline{U_{2m-1}})\subset U_{2m}\,,
\quad \gamma_m^{-1}(C\setminus \overline{U_{2m}})\subset U_{2m-1}\qquad
\mbox{for} \quad  1\le m\le l\,,\qquad \mbox{and}\\
&&\gamma_{m}\big(K/M\setminus U_{m})\subset
U_{m}\,,\quad\gamma_{m}^{-1}\big(K/M\setminus U_{m})\subset
U_{m}\quad\quad\mbox{for}\quad  l+1\le m\le l+p\,.\ee
In particular, the finitely generated group $\langle
\gamma_1,\gamma_2,\ldots,\gamma_l,\gamma_{l+1},\ldots,\gamma_{l+p}\rangle
\subset\is^o(\XX)$ is free and discrete.
\end{thr}
\prf\  For $1\le i\le 2l+p$ we choose open neighbourhoods $U_i\subset K/M$ of $\piF(\xi_i)$ \st
$$ \overline{U_i}   \subset  \bigcap_{\begin{smallmatrix}n=1\\
n\neq i\end{smallmatrix}}^{2l+p} \horF(\xi_n)\,, \quad\mbox{and put}\quad C:=\bigcup_{i=1}^{2l+p} \overline{U_i}\subset K/M\,.$$
Since for any $m\in\{1,2,\ldots,l\}$ we have $\xi_{2m}\in\regrand$ and
$\piF(\xi_{2m-1})\in\horF(\xi_{2m})$,  there exist regular unit speed
geodesics $\sigma_m:\;\RR\to\XX$ \st
$\sigma_m(\infty)=\xi_{2m}$ and $\piF(\sigma_m(-\infty))=\piF(\xi_{2m-1})$. For $1\le m\le l$ let $\gax_m$ be a regular axial isometry with $ \gax_m \sigma_m(t)=\sigma_m(t+1)$ for all $t\in\RR$. In particular, $\gax_m$ possesses the attractive and repulsive fixed points
$\gax_m^+=\sigma_m(\infty)=\xi_{2m}\,$ and
$\gax_m^-=\sigma_m(-\infty)$. By
  Corollary~\ref{dist} and the fact that  $C\setminus
  \overline{U_{2m-1}}\subset\horF(\xi_{2m-1})$, $C\setminus
  \overline{U_{2m}}\subset\horF(\xi_{2m})$ are compact, there exists 
$k_{m}\in\NN$ \st  for all  $j\ge k_m$ we have $\gax_{m}^{j}(C\setminus
  \overline{U_{2m-1}})\subset U_{2m}$
  and $\gax_{m}^{-j}(C\setminus
  \overline{U_{2m}})\subset
  U_{2m-1}$. Putting $\gamma_m:=\gax_m^{k_m}$, we obtain the
  desired regular axial isometries and the corresponding open neighbourhoods.

Similarly, for $1\le m\le p$ we choose a generic
  parabolic isometry $u_m$ in the horospherical subgroup associated to
  $\xi_{l+m}$. By  Proposition~\ref{dyngenuni},
there exists $k_m\in\NN$ \st 
  for all  $j\ge k_m$ we have
$$ u_m^{j}(K/M\setminus U_{l+m})\subset U_{l+m}\qquad\mbox{and}\quad u_m^{- j}(K/M\setminus U_{l+m})\subset U_{l+m}\,.$$ 
Putting $\gamma_{l+m}:=u_m^{k_m}$ we obtain the
  desired generic parabolic isometries with corresponding open neighbourhoods
  $U_{l+1},\ldots, U_{l+p}$.

In order to apply  Klein's Criterion, we put $S_1:=U_1\cup U_2$ and
$S_2:=U_3\cup U_4 $. Since $S_2\subset C\setminus
(\overline{U_1}\cup\overline{U_2})$
we have $\langle \gamma_1\rangle\at S_2\subset  S_1$. Similarly,  $S_1\subset C\setminus
(\overline{U_3}\cup\overline{U_4})$  implies $\langle \gamma_2\rangle\at S_1\subset  S_2$. Hence the group generated by $\gamma_1$ and $\gamma_2$
  is free by Klein's Criterion. 

For $i\in\{1,\ldots,l+p\}$ we denote $\Gamma_i$ the group
  generated by the elements $\gamma_m$ for $m\le i$. If $2\le i\le l-1$
we put $S_i':=\bigcup_{n=1}^i (U_{2n-1}\cup U_{2n})\,$  and 
$S_{i+1}:=
U_{2i+1}\cup U_{2i+2}\,.$
Since $$S_{i+1}\subset C\setminus \bigcup_{n=1}^{2i}\overline{U_{n}}$$
we have  $\gamma\cdot S_{i+1}\subset S_i'$ for all
$\gamma\in\Gamma_i\setminus\{\id\}$. From 
$S_i'\subset C \setminus (\overline{U_{2i+1}}\cup\overline{U_{2i+2}})$
we further obtain $\langle \gamma_{i+1}\rangle \at S_i'\subset S_{i+1}$, hence the group $\Gamma_{i+1}$ generated by the elements $\gamma_m$ for  $m\le i+1$,
  is free. We conclude inductively that $\Gamma_l:=\langle
  \gamma_1,\gamma_2,\ldots,\gamma_l\rangle$ is free.

We next consider the sets $S_l'$ and $U_{2l+1}$. Since $S_l'\subset
 K/M\setminus U_{2l+1}$ and 
$$U_{2l+1}\subset C \setminus \bigcup_{n=1}^{2l}\overline{U_{n}}$$
we have $\langle \gamma_{l+1}\rangle\at S_l'\subset U_{2l+1}$  and $\gamma\cdot U_{2l+1}\subset S_l'$ for all
$\gamma\in\Gamma_l\setminus\{\id\}$, hence the group $\Gamma_{l+1}$ is free. 
 
For $m\in\{2,\ldots,p-1\}$ we put $ S_{l+m}':=\bigcup_{n=1}^{2l+m} U_n\,.$
Again from  $ S_{l+m}'\subset K/M\setminus U_{2l+m+1}$ and
$$ U_{2l+m+1}\subset C \setminus \bigcup_{n=1}^{2l+m}\overline{U_{n}}$$
we obtain $\langle \gamma_{l+m+1}\rangle\at S_{l+m}'\subset U_{2l+m+1}$ and $\gamma\cdot U_{2l+m+1}\subset S_{l+m}'$ for all
$\gamma\in\Gamma_{l+m}\setminus\{\id\}$.  We conclude inductively that $\langle
  \gamma_1,\gamma_2,\ldots,\gamma_{l+p}\rangle$ is free.

Finally suppose $\Gamma:=\langle
  \gamma_1,\gamma_2,\ldots,\gamma_{l+p}\rangle$ is not discrete. Then there exists
  a sequence $(\gax_j)\subset\Gamma$
  converging to the identity. For $j\in \NN$ we write $\gax_j:=
  s_{k_j}^{(j)}s_{k_j-1}^{(j)}\ldots s_{1}^{(j)}$ as a reduced word,
  i.e. $s_m^{(j)}\in \{
  \gamma_1,\gamma_1^{-1},\gamma_2,\gamma_2^{-1}\ldots,\gamma_{l+p},\gamma_{l+p}^{-1}\}$
  for $1\le m\le k_j$, and $s_{m+1}^{(j)}\ne
  (s_{m}^{(j)})^{-1}$ for $1\le m\le k_j-1$.  Passing to a subsequence if
  necessary, we may assume that $s_{k_j}^{(j)}=s\,$ and $s_{1}^{(j)}=s'$ for all
  $j\in\NN$. We denote $U, U'\in
  \{U_1, U_2,\ldots, U_{2l},\ldots, U_{2l+p}\}$ the corresponding 
  neighbourhoods of the (attractive) fixed point of $s$, $s'$.
  Let $\eta\in\regrand$ \st 
$\pi^B(\eta)\in C\setminus (\overline{U}\cup \overline{U'})$.
Then  $\piF(s_1^{(j)}\eta)\in U'$, 
hence by the dynamics of the generators of $\Gamma$ the point $\pi^B(\gax_j\eta)$ is contained in 
  $U$ for all $j\in\NN$. On the other hand, since $\gax_j$ converges to the identity, $\pi^B(\gax_j\eta)$ converges to
  $\pi^B(\eta)\in C\setminus (\overline{U}\cup \overline{U'})$,
 a contradiction.  \qed

\section{Discrete isometry groups}\label{GeomLimit}

In this section we define the geometric limit set of a discrete isometry group
of a globally symmetric space $\XX$ of noncompact type. We extend the familiar
notion of ``nonelementary'' groups from rank one to higher rank symmetric spaces.  We will then describe the structure of the limit set of such nonelementary groups using {\rm Theorem}~\ref{hypdyn} as our main tool. 

In this section, $\XX$ will again denote a globally symmetric space of
noncompact type with base point $\xo\in\XX$, and $G=\is^o(\XX)$ the connected
component of the identity.

\subsection{The limit set}

\begin{df} 
A subgroup $\Gamma\subset G$ is called {\hd discrete}, if it has a discrete orbit in $\XX$. In this case, the {\hd geometric limit set} $L_\Gamma$ of $\Gamma$ is
defined by 
$\  \Lim :=\Gamma\at\xo \cap \rand$. 
\end{df}

We remark that this definition can be extended to isometry groups of arbitrary
Hadamard manifolds. Furthermore, the geometric limit set does not depend on the chosen base point $\xo$. 

Due to the rich algebraic structure of symmetric spaces, we can consider
various sets describing asymptotic properties of a discrete isometry group
$\Gamma$. In order to do so, we fix a Cartan decomposition   \Cart\ \wrt
$\xo\in\XX$. Recall that the map $\pi^B:\,\regrand\to K/M\ $ denotes the projection introduced in
section~\ref{compactify}. 
\begin{df}\label{varlimsets}
We call the set of regular limit points $\reglim:=\Lim\cap\regrand$ the {\hd
  regular limit set}, and the projection $K_\Gamma:=\pi^B(\reglim)\subseteq
K/M$ the {\hd transversal limit set} of $\Gamma$. The  set of Cartan
projections of all points in the geometric limit set $\Lim$ of $\Gamma$ is
called the {\hd directional limit set} $P_\Gamma\subseteq\overline{\aL^+_1}$,
the {\hd limit cone} $\ell_\Gamma\subseteq\aL^+_1$ of $\Gamma$ is defined by
$$ \ell_\Gamma:=\{ L(\gamma)/||L(\gamma)||\;|\, \gamma\in \Gamma\
\mbox{regular axial}\}\,. $$
\end{df} 

The following lemma is well-known and remains true for discrete isometry groups of Hadamard manifolds. We include the proof for the convenience of the reader.
\begin{lem}\label{discretecommute}
Let $\Gamma\subset G$ be a discrete group and $\gax\in\Gamma$ axial. If $\varphi\in\Gamma$ fixes $\gax^+$, then $\varphi$ commutes with a power of $\gax$. 
\end{lem}
\prf\ Let $x\in\Ax(\gax)$ and $l>0$ the translation length of $\gax$.  Then for all $n\in\NN$ 
$$d(x,\gax^{-n}\varphi\gax^n x)= d(\gax^n x, \varphi \gax^n x)  
= d(\sigma_{x,\gax^+}(n l),\varphi\sigma_{x,\gax^+}(nl))$$ 
is bounded from above by a constant $r\ge 0$ because $\varphi$ fixes $\gax^+$. Since 
$ \# (\Gamma x\cap B_x(r))$ is finite by discreteness of $\Gamma$, we conclude that
$$\gax^{-n}\varphi \gax^n=\gax^{-m}\varphi \gax^m $$ for integers $m>n$, i.e. $\varphi$ commutes with $\gax^{m-n}\neq\id$.\qed

\subsection{Nonelementary groups}\label{nonel}

We are now going to generalise to symmetric spaces $\XX$ of higher rank
the notion of ``nonelementary groups'' familiar in the context of isometry
groups of real hyperbolic spaces. Let \Cart\ be a Cartan decomposition \wrt $\xo\in\XX$.
\begin{df}
A discrete subgroup $\Gamma$ of the isometry group $\is^o(\XX)$ is
called {\hd nonelementary} if  $\reglim\ne \emptyset$ and if for all 
$ \xi$, $\eta\in \reglim$  we have
$$\pi^B(\Gamma\at \xi)\cap \horF(\eta)\neq \emptyset\,. $$
Otherwise $\Gamma$ is called {\rm elementary}.
\end{df}
Notice that an abelian discrete group $\Gamma\subset\is^o(\XX)$ of axial
isometries is 
elementary, because its limit set is contained in the boundary of the
invariant maximal flats. Hence $\Gamma\at\xi=\xi$ for every $\xi\in\reglim$ which
implies $\pi^B(\Gamma\at \xi)=\pi^B(\xi)\notin \horF(\xi)$. The same argument
shows that a discrete group $\Gamma\subset\is^o(\XX)$ is elementary, if it is
contained in the  stabiliser of a regular limit point. 

The following lemma provides a first example of nonelementary groups:
\begin{lem}\label{freeisnonel}
If $\Gamma\subset G$ is a discrete group containing $l\ge 2$ regular axial
isometries $\gamma_1,\gamma_2,\ldots,\gamma_l$ with 
\begin{eqnarray}\label{addproperty}
 \pi^B(\gamma_i^+),\,\pi^B(\gamma_i^-) & \in & 
  \bigcap_{\begin{smallmatrix}m=1\\
m\neq i\end{smallmatrix}}^l\left( \horF(\gamma_m^+)\cap
  \horF(\gamma_m^-)\right)\qquad\mbox{for}\quad 
  1\le i\le l\,, \quad \mbox{and}\nonumber\\
K_\Gamma & \subseteq & \bigcup_{i=1}^l \big(\horF(\gamma_i^+)\cup
\horF(\gamma_i^-)\big)\,,
\end{eqnarray} 
then $\Gamma$ is nonelementary.

In particular, free groups generated by $l\ge 2$ regular axial isometries $\gamma_{1},\ldots,\gamma_{l}$ as in Theorem~\ref{freeconst} and the additional property (\ref{addproperty}) are nonelementary.
\end{lem}
\prf\  We denote $\Gamma':=\langle \gamma_1,\gamma_2,\ldots,\gamma_l\rangle$
and   let 
$U_1,U_2,\ldots, U_{2l}\subset K/M$ be pairwise disjoint
open sets as in {\rm Theorem}~\ref{freeconst}.  Choose $\xi$,
$\eta\in\reglim$ arbitrary. Then by (\ref{addproperty}) there exist $g,h \in\{\gamma_1,\gamma_1^{-1},\ldots,\gamma_l,\gamma_l^{-1}\}$ \st
$\piF(\xi)\in\horF(g^+)$ and $\piF(\eta)\in\horF(h^+)$.

If $g\ne h$, then $\piF(g^-)\in\horF(h^-)$. Since $\horF(h^-)$ is open, and
$\pi^B(g^{-j}\xi)$ converges to $\pi^B(g^-)$ as $j\to\infty$ by {\rm Corollary}~\ref{dist}, there exists
$n\in\NN$ \st $\pi^B(g^{-n}\xi)\in \horF(h^-)$. Now $\piF(h^+)\in \horF(\eta)$
and the same argument as before imply the existence of $m\in\NN$ \st
$\pi^B(h^{m}g^{-n}\xi)\in\horF(\eta)$.

If $g=h$, we choose
$\gamma\in\{\gamma_1,\gamma_1^{-1},\ldots,\gamma_l,\gamma_l^{-1}\}\setminus\{g,g^{-1}\}$.
Then we have $\piF(\xi)\in\horF(g^+)$, $\piF(g^-)\in\horF(\gamma^-)$, 
$\piF(\gamma^+) \in\horF(g^-)$ and $\piF(g^+)\in\horF(\eta)$. Therefore {\rm
  Corollary}~\ref{dist} implies the existence of integers $n,m,k\in\NN$ \st
$\piF(g^{-n}\xi)\in\horF(\gamma^-)$, $\piF(\gamma^{m} g^{-n}\xi)\in\horF(g^-)$ and
finally $\piF(g^{k}\gamma^{m} g^{-n}\xi)\in\horF(\eta)$. \qed \\

We remark that condition~(\ref{addproperty}) is not very restrictive, because the sets $\horF(\gamma_i^+)$ are dense and open in $K/M$.  
Moreover, the same arguments as in the proof above yield the following lemma:
\begin{lem}
Let $\Gamma\subset G$ be a discrete group containing a generic parabolic isometry $u$ with fixed point $\eta\in\regrand$. If $\Gamma$ further contains a regular axial isometry $\gax$ \st $ \piF(\eta)\in \horF(\gax^+)\cap\horF(\gax^-)$, or a second generic parabolic isometry with fixed point $\zeta\in\regrand$ and 
$\piF(\zeta)\in\horF(\eta)\,,$
then $\Gamma$ is nonelementary. 

In particular, free groups $\langle \gamma_1,\ldots,\gamma_l,\gamma_{l+1},\ldots,\gamma_{l+p}\rangle$ as in Theorem~\ref{freeconst} with at least two generators and $p\ge 1$ are nonelementary. 
\end{lem}

We are finally able to justify our choice of the notion ``nonelementary''.
The following lemmata state that for torsion free discrete isometry groups of rank one symmetric spaces, our definition coincides with the familiar one.  
\begin{lem}\label{classic}
If $\rank(\XX)=1$, then a discrete isometry group
  $\Gamma\subset\is^o(\XX)$ is nonelementary 
  if it possesses infinitely many limit points.
\end{lem}
\prf\  Since $\rank(\XX)=1$, we have $\regrand=\rand$ by convention, $\rand$ is homeomorphic to $\pi^B(\rand)$, and $\rand=\horinf(\zeta)\cup
  \{\zeta\}$ for any point $\zeta$ in the
  geometric boundary. 

Suppose $\Gamma\subset G=\is^o(\XX)$ possesses infinitely many limit
  points, and assume there exist $\xi,\eta\in L_{\Gamma}$ \st 
  $\Gamma\at\xi\cap\horinf(\eta)=\emptyset$. Then $\gamma\xi=\eta$ for
  all $\gamma\in\Gamma$, in particular $\xi=\eta$. This implies that every
  element in $\Gamma$ fixes $\xi$. Let $\Gamma'\subseteq\Gamma$ be a torsion
  free subgroup of finite index which exists by Selberg's Lemma. Since
  $\Gamma'$ does not contain elliptic elements, $\Gamma'$ contains only
  parabolic and axial isometries which all fix $\xi$. By Lemma~\ref{discretecommute}, the set
  of axial elements in $\Gamma'$ must all have the same axis. We conclude that
  $\Gamma'$ possesses at most two limit points, hence $\Gamma$ possesses only
  finitely many limit points, a  contradiction.\qed

\begin{lem}\label{classic}
If $\rank(\XX)=1$, then a torsion free nonelementary discrete isometry group
  $\Gamma\subset\is^o(\XX)$ possesses infinitely many limit points.
\end{lem}
\prf\ Suppose  $\Gamma$ possesses only finitely many limit points. Since
$\reglim\ne \emptyset$ by definition, we first treat the case
$\Lim=\{\xi\}$. But then the $\Gamma$-invariance of the limit set implies
$\Gamma\at \xi=\xi$, hence in particular $\Gamma\at\xi\cap\horinf(\xi)=\emptyset$, a
contradiction.

Next suppose that $\Lim=\{\xi,\eta\}$ with $\eta\ne\xi$. By $\Gamma$-invariance of the limit
set we conclude that $\Gamma$ leaves invariant the unique geodesic $\sigma$
joining $\xi$ to $\eta$. Since $\Gamma$ does not contain elliptic elements, we
conclude that $\Gamma$ is a cyclic group generated by an axial element which
translates $\sigma$, hence $\Gamma\at\xi=\xi$ and $\Gamma\at\eta=\eta$. Again
$\Gamma\at\xi\cap\horinf(\xi)=\emptyset$ gives a contradiction. 

Finally, if $\Lim$ contains at least three points $\xi,\eta,\zeta$, then $\Gamma$ contains either
three parabolics $g,\gamma,p$ \st $g\xi=\xi$, $\gamma\eta=\eta$ and $p\zeta=\zeta$, or
an axial isometry $\gax$ and a second isometry $\gamma$ which does not fix
$\gax^+$ or $\gax^-$. In the first case,
the dynamics of parabolic isometries of rank one symmetric spaces imply that
the points $p^j\xi$ for $j\in\ZZ$ are all disjoint, and contained in the limit
set. Hence $\#\Lim=\infty$. 

In the second case, if $\xi\in\Lim$ denotes the point not fixed by $\gax$, the
points $\gax^j\xi$, $j\in\ZZ$, are disjoint and  belong to the limit set,
hence again $\#\Lim=\infty$.\qed \\

We remark that the restriction to torsion free groups is only necessary for
the following particular situation which may occur: If $\Gamma$ consists of an
axial isometry $\gax$ which translates a geodesic $\sigma$, and an elliptic
isometry $e$ which leaves invariant  $\sigma$ and permutes the extremities of
$\sigma$, then $\langle \gax, e\rangle$ is nonelementary \wrt our
definition, but possesses only two limit points. 

The above proof and the remark before Lemma~\ref{freeisnonel} show that a
torsion free discrete isometry group of a rank one symmetric space is
elementary if and only if it is contained in the stabiliser of a limit point. 

\subsection{The approximation argument}

Since we do not know much about the dynamics of parabolic isometries, the
description of the limit set of discrete isometry groups which, in general,
contain parabolics, is difficult. Fortunately, the following ``approximation argument'' allows to approach every regular limit point in a nonelementary discrete group by a sequence of regular axial isometries. The main tool in the proof is Proposition~4.5.14 in \cite{E}. 

\begin{prp}\label{axialapprox}
Let $\Gamma\subset G =\is^o(\XX)$ be a nonelementary discrete group. Then for every
$\xi\in\reglim$ there exists a sequence of axial isometries
$(\gax_j)\subset\Gamma$ 
\st $\gax_j^+$ converges to $\xi$ and $\gax_j^{-}$ converges to a
point in $\horinf(\xi)$.  Furthermore, $d(\xo,\Ax(\gax_j))$ is bounded as $j\to\infty$.
\end{prp}
\prf\  Fix $\xi\in\reglim$ and let $\zeta\in \reglim$ arbitrary. Since $\Gamma$ is nonelementary, there exists $\gamma\in\Gamma$ \st $\pi^B(\gamma\zeta)\in\horF(\xi)$. This implies the existence of a regular unit speed geodesic $\sigma$ with $\sigma(\infty)=\xi$ and $\pi^B(\sigma(-\infty))=\pi^B(\gamma\zeta)$. 

Let $\Phi_t:\RR\to S\XX$ denote the geodesic flow of $\XX$ and $p:S\XX\to\XX$ the foot point projection. By Proposition~4.5.14 in \cite{E}, there exists a sequence of unit tangent vectors $(v_j)\subset S\XX$ converging to the tangent vector $\dot\sigma(0)\in S\XX$ of $\sigma(0)$, and a sequence $(\gax_j)\subset\Gamma$ of regular axial isometries with translation lengths $(l_j)\subset\RR$ \st 
$$ \gax_j \Phi_tv_j=\Phi_{t+l_j}v_j\qquad \mbox{for all}\quad j\in\NN\,.$$
Hence $\gax_j^+\to\sigma(\infty)=\xi$, $\gax^{-}_j\to\sigma(-\infty)$, and
$$d(\xo,\Ax(\gax_j))\le d(\xo, \sigma(0))+d(\sigma(0),p(v_j))\to d(\xo,\sigma(0)) \qquad\mbox{as}\quad j\to \infty\,.$$ 
Hence $(\gax_j)$ is the desired sequence. \qed\\

An immediate corollary of this proposition is the following
\vspace{-2mm}\begin{thr}\label{axdens}
If $\Gamma\subset G=\is^o(\XX)$ is a nonelementary discrete group, then the set of attractive fixed
points of regular axial isometries is a dense subset of the limit set $\Lim$.
\end{thr}

\subsection{Sequences of axial isometries}

From here on $\Gamma\subset G$ will always denote a nonelementary discrete
group. The following equivalence for certain sequences of axial isometries in
$\Gamma$ will be necessary in the proof of Theorem~\ref{hypdyn}.
\begin{lem}\label{axeq}
Fix $x\in\XX$ and let $(\gax_j)\subset\Gamma$ be a sequence of regular axial isometries \st
 $d(x,\Ax(\gax_j))$ remains bounded as $j\to\infty$. Then $(\gax_j x)\subset\XX$ converges to a boundary point $\xi\in\rand$ in the cone
 topology if and only if the 
 sequence of attractive fixed
 points $(\gax_j^+)\subset\regrand$ of $(\gax_j)$ converges to $\xi$ in the cone topology. 
\end{lem}
\prf\  Let $(\gax_j)$ be a sequence of axial isometries with attractive fixed
points $(\gax_j^+)\subset\regrand$, and $c\ge 0$ \st
$d(x,\Ax(\gax_j))\le c$. For all $j\in \NN$ we let $x_j\in\Ax(\gax_j)$ be the orthogonal projection of $x$ to $\Ax(\gax_j)$, and $l_j:=d(x_j,\gax_j x_j)$ the translation length of $\gax_j$.

It suffices to prove that for all $R>>1$
$$d(\sigma_{x,\gax_j x}(R), \sigma_{x,\gax_j^+}(R))\to 0\qquad\mbox{as}\quad j\to\infty\,.$$
For $j\in\NN$ we put $d_j:=d(x,\gax_j x)$. Using the convexity of the distance function, the triangle inequality and $|d_j-l_j|\le 2c$,  we compute
\be&& d(\sigma_{x,\gax_jx}(R), \sigma_{x,\gax_j^+}(R))\\
&&\quad \le\  \frac{R}{l_j}\big(d(\sigma_{x,\gax_j x}(l_j),\sigma_{x,\gax_jx}(d_j))+d(\sigma_{x,\gax_jx}(d_j),\gax_j x _j) + d(\gax_j x_j,\sigma_{x,\gax_j^+}(l_j))\big)\\
&& \quad =\  
\frac{R}{l_j}\big(|l_j-d_j| +d(\gax_jx,\gax_j x_j)
\  + d(\sigma_{x_j,\gax_j^+}(l_j),\sigma_{x,\gax_j^+}(l_j))\big) \le 4c\frac{R}{l_j}\,.
\ee
From  $l_j\to\infty$ as $j\to\infty$ we conclude that $d(\sigma_{x,\gax_j x}(R), \sigma_{x,\gax_j^+}(R))$ tends to zero.\qed

For the remainder of this section we fix a Cartan decomposition
\Cart\ \wrt  $\xo\in\XX$. The following result generalises Corollary~\ref{dist} to sequences of axial isometries $(\gax_j)\subset \Gamma$. In combination with the approximation argument Proposition~\ref{axialapprox} this will be our main tool for the description of the structure of the limit set in the following section. 
\begin{thr}\label{hypdyn} 
Let  $(\gax_j)\subset\Gamma$ be a sequence of regular axial isometries \st $\gax_j\xo$
 converges to $\xi^+\in\regrand$, $\gax_j^{-1}\xo$ converges
 to  $\xi^-\in\horinf(\xi^+)$, and 
 $d(\xo,\Ax(\gax_j))$ is bounded as $j\to\infty$.  Then for 
 any $\zeta\in\regrand$ with $\piF(\zeta)\in\horF(\xi^-)$ there exist integers $n_j$, $j\in\NN$, \st the sequence 
$(\gax_j^{n_j} \zeta)$ converges to the unique point $\eta^+\in G\at\zeta$
 with $\piF(\eta^+)=\piF(\xi^+)\,.$ 
In particular, if $\zeta\in\horinf(\xi^-)$, then there exist integers $n_j$,
 $j\in\NN$, \st $\gax_j^{n_j} \zeta$
 converges to $\xi^+$.
\end{thr}  
\prf\   Let $(\gax_j)$ be a sequence of regular axial isometries with the
properties stated in the theorem. We denote by $\gax_j^+$ the
attractive fixed point of $\gax_j$ and by $H_j\in\aL^+_1$ its Cartan projection. Let $\zeta\in\regrand$ with
$\piF(\zeta)\in\horF(\xi^-)$ arbitrary, and $\eta^+\in G\at\zeta$ the unique
point \st $\piF(\eta^+)=\piF(\xi^+)$. By the previous lemma,  $\gax_j^+$
converges to $\xi^+$ in the cone topology, hence Lemma~\ref{coneKM} implies
that $\pi^B(\gax_j^+)$
converges to $\piF(\xi^+)$ in $K/M$, and the Cartan projections
$(H_j)\subset\aL^+_1$ converge to the Cartan projection $H_\xi$ of $\xi^+$. 

Let $Y\subset G\at\zeta$ be an open neighbourhood of $\eta^+$. Then for $j>N$
sufficiently large, we have $\piF(\gax_j^+)\in\piF(Y)$ and
$\piF(\zeta)\in\horF(\gax_j^-)$, hence by {\rm Corollary}~\ref{dist} there
exists $n_j\in\NN$ \st
$\gax_j^{n_j}\zeta\in Y$. We conclude that the sequence $(\gax_j^{n_j}\zeta)$ converges to $\eta^+$ in the cone
topology. \qed\\

\subsection{The structure of the limit set}\label{Limsetstructure}

We are finally able to describe precisely 
the limit set of nonelementary discrete groups acting on a globally symmetric
space $\XX$ of noncompact type. We fix a Cartan decomposition \Cart\ \wrt 
a base point $\xo\in\XX$. 
Recall the definitions of $K_\Gamma$, $P_\Gamma$ and $\ell_\Gamma$ from  Definition~\ref{varlimsets}. 
\begin{thr}\label{minclosed}
Let $\Gamma\subset G= \is^o(\XX)$ be a nonelementary discrete group. Then the transversal
limit set $K_\Gamma$ is a minimal closed set under the action of
$\Gamma$. 
\end{thr}
\prf\  Since $\Gamma$ is nonelementary, $K_\Gamma$ is nonempty. 
We fix $k_0 M\in K_\Gamma$, and let $k M \in K_\Gamma$ 
arbitrary. Let $\xi^+\in\reglim$ \st $\pi^B(\xi^+)=kM$ and
denote by $H\in\aL^+_1$ the Cartan projection of $\xi^+$. Due to
{\rm Proposition}~\ref{axialapprox} there exists a sequence
$(\gax_j)\subset \Gamma$ of axial isometries \st $\gax_j^+$ converges to
$\xi^+$ and $\gax_j^{-}$ converges to a point
$\xi^-\in\horinf(\xi^+)$. Put $\xi_0:=(k_0,H)$ with
$H\in\aL^+_1$ as above, i.e. $\xi_0\in G\at\xi^+$. If $k_0M\in\horF(\xi^-)$, then $\xi_0\in\horinf(\xi^-)$ and, by {\rm
  Theorem}~\ref{hypdyn}, there exist integers $n_j$, $j\in\NN$, \st
$$ \lim_{j\to\infty} \gax_j^{n_j} \xi_0 =\xi^+ \,.$$
If $k_0 M =\pi^B(\xi_0)\notin \horF(\xi^-)$, there exists $\gamma\in\Gamma$  \st $\pi^B(\gamma\xi_0)\in \horF(\xi^-)$ because $\Gamma$ is nonelementary. Therefore
$\gamma\xi_0\in\horinf(\xi^-)$ and {\rm
  Theorem}~\ref{hypdyn} guarantees the existence of integers $n_j$, $j\in\NN$, \st
$$ \lim_{j\to\infty} \gax_j^{n_j} (\gamma \xi_0)=\xi^+\,.$$
This shows that
$\Gamma(k_0 M)=\pi^B(\Gamma\xi_0)$ is dense in $K_\Gamma$, hence its closure
is a minimal closed set under the action of $\Gamma$. \qed\\

For nonelementary Schottky groups as in section~\ref{free} this implies the
following
\begin{cor}
If $\Gamma:=\langle
\gamma_1,\ldots,\gamma_l,\gamma_{l+1},\ldots,\gamma_{l+p}\rangle$ is a
nonelementary Schottky group with corresponding open neighbourhoods
$U_1,U_2,\ldots, U_{2l},U_{2l+1},\ldots, U_{2l+p}\subset K/M$ as in
Theorem~\ref{freeconst}, then 
$$K_\Gamma\subset \bigcup_{i=1}^{2l+p} U_i\,.$$
\end{cor}
\prf\  If $l=0$, we let $\xi\in\regrand$ be a fixed point of the generic
  parabolic isometry $\gamma_{l+1}=\gamma_1$, if $l\ge 1$, we put
  $\xi:=\gamma_1^+$. Then $\xi\in\reglim$, and by Theorem~\ref{minclosed} we
  have $ K_\Gamma= \piF(\overline{\Gamma\at\xi})$. Let
  $(\gax_j)\subset\Gamma$ be an arbitrary sequence. For $j\in\NN$ we write $\gax_j:=
  s_{k_j}^{(j)}s_{k_j-1}^{(j)}\ldots s_{1}^{(j)}$  as a reduced word,
  i.e. $s_m^{(j)}\in \{
  \gamma_1,\gamma_1^{-1},\gamma_2,\gamma_2^{-1}\ldots,\gamma_{l+p},\gamma_{l+p}^{-1}\}$
  for $1\le m\le k_j$, and $s_{m+1}^{(j)}\ne
  (s_{m}^{(j)})^{-1}$ for $1\le m\le k_j-1$.  Passing to a subsequence if
  necessary, we may assume that $s_{k_j}^{(j)}=s\,$ for all
  $j\in\NN$. Let $U\in
  \{U_1, U_2,\ldots, U_{2l},\ldots, U_{2l+p}\}$ be the corresponding 
  neighbourhood of the (attractive) fixed point of $s$. Notice that if
  $s_1^{(j)}\in\{\gamma_1,\gamma_1^{-1}\}$, then $s_1^{(j)}\xi=\xi$. Hence for
  $j\in\NN$ we
  put $i_j:=\min\{i\in\NN\;|\, s_i^{(j)}\notin\{\gamma_1,\gamma_1^{-1}\}\}$.
 Then $\gax_j\xi=s\at s_{k_j-1}^{(j)}\ldots s_{i_j}^{(j)}\xi$, and from the
  dynamics of the generators of $\Gamma$ we have
  $\pi^B(\gax_j\xi)\in U$ for all $j\in\NN$. We conclude that
  $\piF(\overline{\Gamma\at\xi})\subset \bigcup_{i=1}^{2l+p} U_i\,.$\qed

\begin{thr}\label{Product}
Let $\Gamma\subset G= \is^o(\XX)$ be a  nonelementary discrete group of
isometries. Then the regular geometric 
limit set is isomorphic to the product $K_\Gamma\times
(P_\Gamma\cap \aL^+_1)$.
\end{thr}
\prf\  If $\xi\in \reglim$, then $\pi^B(\xi)\in K_\Gamma$ and the Cartan projection of $\xi$ belongs to $P_\Gamma\cap\aL^+_1$. 

Conversely, let $kM\in K_\Gamma$ and $H\in P_\Gamma\cap\aL^+_1$. By definition
of $P_\Gamma$, there
exists a sequence $(\gamma_j)\subset\Gamma$  \st the Cartan projections $(H_j)\subset\overline{\aL^+}$ of
$(\gamma_j\xo)$ satisfy $\angle (H_j, H)\to 0$ as $j\to
\infty$.
Furthermore, a subsequence of $(\gamma_j\xo)$ converges to a point  $\xi_0=(k_0,H)\in\reglim$ where $k_0\in K$. 

By {\rm Theorem}~\ref{minclosed}, $K_\Gamma=\overline{\Gamma(k_0M)}$ is a
minimal closed set under the action of $\Gamma$, hence 
$$kM\in\overline{\Gamma(k_0M)}=\pi^B(\overline{\Gamma\at\xi_0})\,.$$
Since the action of $\is^o(\XX)$ on the geometric boundary does not change
 Cartan projections, we conclude 
that  the closure of $\Gamma\at \xi_0$ contains $(k,H)$. In particular,
$(k,H)\in\overline{\Gamma\at\xi_0}\subseteq\reglim$.\qed\\

The limit cone and the directional
limit set of a discrete isometry group are related as follows:
\begin{thr}\label{Cone}
If $\Gamma\subset G= \is^o(\XX)$ is a nonelementary discrete group, then 
 $P_\Gamma=\overline{\ell_\Gamma}$.
\end{thr}
\prf\  In order to prove $P_\Gamma\supseteq\overline{\ell_\Gamma}$, we let
$H\in\overline{\ell_\Gamma}$ arbitrary.  Then there
exists a sequence $(\gax_j)\subset\Gamma$ of axial isometries with
translation directions $L_j:=L(\gax_j)$ satisfying $\angle L_j,H)\to
0$ as $j\to\infty$. 

Suppose $H\notin P_\Gamma$, and, for $\gamma\in\Gamma$, let $H_\gamma$ denote the Cartan projection of
$\gamma\xo$. Then there exists $\eps>0$  \st $\angle
(H_\gamma,H)>\eps\,$ for all but finitely many $\gamma\in\Gamma$.
Put $\xi=(\id,H)\in\rand$. For $j\in\NN$, we let 
$g_j\in G$ \st $\gax_j=g_j e^{L_j}g_j^{-1}$, put $x_j:=g_j \xo\in
 \Ax(\gax_j)$, and let $n_j$ be an integer greater than $2j
 d(\xo,x_j)/||L_j||$. We abbreviate $\gamma_j:=\gax_j^{n_j}=g_j e^{L_j n_j}g_j^{-1}$ and 
 $H_j:=H_{\gamma_j}$. By 
$G$-invariance of the directional distance and equation~(\ref{dirwinkel}) we have
\be \bs_{G\cdot\xi}(x_j,\gamma_j x_j)&=&\bs_{G\cdot\xi}(g_j^{-1} x_j,g_j^{-1} \gamma_j
  x_j)= \bs_{G\cdot\xi}(\xo,e^{L_j  n_j}\xo)\\
&=& \langle H, L_j  n_j\rangle =  n_j \langle  L_j,H \rangle\,,\ee
and $\angle(L_j,H)\to 0$ implies $\langle L_j/\Vert L_j\Vert, H\rangle \to 1$
as $j\to\infty$. 
Using again equation~(\ref{dirwinkel}), the $G$-invariance   and the
triangle inequality  for the Riemannian and the directional distance, and
$2 d(\xo,x_j)\le n_j\Vert L_j\Vert/j$
 we conclude
\be \cos \angle (H_j, H)&=&\frac{\langle H_j,H\rangle}{\Vert H_j\Vert }
=  \frac{\bs_{G\cdot\xi}(\xo,\gamma_j \xo)}{d(\xo,\gamma_j\xo)}
\ge   \frac{\bs_{G\cdot\xi}(x_j,\gamma_j
 x_j)-2d(\xo,x_j)}{d(x_j,\gamma_j x_j)+2 d(\xo,x_j)}\\
&\ge &\frac{n_j\langle L_j, H\rangle - n_j\Vert L_j\Vert/j}{n_j\Vert L_j\Vert
  +n_j\Vert L_j\Vert /j}=\frac{\langle L_j/\Vert L_j\Vert, H\rangle -1/j}{1+1/j}\ \to\  1\ee
as $j\to\infty$, a contradiction to our assumption.

Conversely, we first prove $P_\Gamma\cap\aL^+_1\subseteq
{\overline\ell_\Gamma}$. 
Given  $H\in P_\Gamma\cap\aL^+_1$, there exists $\xi\in\reglim$ with Cartan
projection $H$. Let $(\gax_j)\subset\Gamma$  be a sequence of regular axial
isometries as in Proposition~\ref{axialapprox} with the
properties $\gax_j^+\to \xi$, $\gax_j^{-}\to\eta\in
\horinf(\xi)$  and $d(\xo,\Ax(\gax_j))\le c$ for some constant $c>0$. Then
Lemma~\ref{axeq} implies $\gax_j\xo\to \xi$, hence by Lemma~\ref{coneKM} the
Cartan projections $(H_j)\subset\aL^+$ of $\gax_j\xo$ satisfy
$$ \langle H_j/\Vert H_j\Vert,H\rangle \to 1\quad \mbox{as}\  j\to\infty\,.$$ 
For $j\in\NN$, we let $x_j\in\Ax(\gax_j)$ be the orthogonal projection of
$\xo$ to $\Ax(\gax_j)$.  If $L_j:=L(\gax_j)\in\overline{\aL^+}$ denotes the
translation vector of $\gax_j$, we estimate
\be \langle \frac{L_j}{||L_j||}, H\rangle&=& \frac{\bs_{G\cdot\xi}(x_j, \gax_j
  x_j)}{d(x_j,\gax_j x_j)} \le\frac{\bs_{G\cdot\xi}(\xo,\gax_j \xo)+2d(\xo,x_j)}{d(\xo, \gax_j \xo)-2d(\xo,x_j)}\\
&\le &  \frac{\langle H_j/\Vert H_j\Vert, H\rangle +2c/\Vert
  H_j\Vert}{1-2c/\Vert H_j\Vert}\
\to\ 1\quad\mbox{and}\\
\langle \frac{L_j}{||L_j||},  H\rangle&\ge&\frac{\bs_{G\cdot\xi}(\xo,\gax_j \xo)-2 d(\xo,x_j)}{d(\xo, \gax_j \xo)}
\ge \langle \frac{H_j}{\Vert H_j\Vert},H\rangle - \frac{2c}{\Vert H_j\Vert}\ \to \ 1\,
\end{eqnarray*}
as $j\to\infty$. This gives  $\angle (L_j, H)\to 0$ as $j\to\infty$, hence
$H\in \overline{\ell_\Gamma}$.  

Since the closure of $P_\Gamma\cap\aL^+_1$ contains $P_\Gamma$, and
${\overline\ell_\Gamma}$ is a closed set in $\overline {\aL^+_1}$, we
conclude that 
$\  P_\Gamma\subseteq\overline{P_\Gamma\cap\aL^+_1}\subseteq{\overline\ell_\Gamma}\,.$ \qed

\vspace{1cm}
Gabriele Link\\
Mathematisches Institut II\\
Universit{\"a}t Karlsruhe\\
Englerstr.~2\\
76 128 Karlsruhe\\
e-mail:\ gabriele.link@math.uni-karlsruhe.de

\end{document}